\documentclass[twoside,leqno,twocolumn]{article}

\usepackage[margin=1in]{geometry}
\usepackage{ltexpprt}
\usepackage{multirow}
\usepackage[ntheorem]{empheq}
\usepackage{verbatim}
\usepackage{graphicx,amsmath,amsfonts,amssymb,subfigure}
\usepackage{tikz}
\usetikzlibrary{arrows,patterns}
\usepackage{mathtools}
\usepackage{booktabs}
\usepackage{textcomp}
\usepackage{setspace}
\usepackage{enumerate}
\usepackage{enumitem}
\usepackage[hidelinks]{hyperref}
\usepackage{pgfplots}
\usepackage{listings}
\usepackage{color} 
\definecolor{mylilas}{RGB}{170,55,241}
\definecolor{myred2}{RGB}{236,1,59}
\definecolor{myblue2}{RGB}{116,122,255}
\definecolor{mygreen2}{RGB}{110,184,129}
\definecolor{mygreen3}{RGB}{80,229,109}
\definecolor{mygray2}{RGB}{196,196,196}
\definecolor{mygreen}{RGB}{28,142,0} 
\definecolor{myc1}{RGB}{0,115,189}
\definecolor{myc2}{RGB}{217,84,26}
\definecolor{myc3}{RGB}{237,176,33}
\definecolor{myc4}{RGB}{59,138,36}
\usetikzlibrary{decorations.markings}
\usepackage[customcolors]{hf-tikz}

\tikzset{style green/.style={
    set fill color=green!50!lime!60,
    set border color=white,
  },
  style cyan/.style={
    set fill color=cyan!90!blue!60,
    set border color=white,
  },
  style mygrey/.style={
    set fill color=white!89.803921568627459!black,
    set border color=white,
  },
  hor/.style={
    above left offset={-0.15,0.31},
    below right offset={0.15,-0.125},
    #1
  },
  ver/.style={
    above left offset={-0.1,0.3},
    below right offset={0.15,-0.15},
    #1
  }
}

\usepackage{colortbl}%
\newcommand{\myrowcolour}{\rowcolor{white!89.803921568627459!black}}

\usepackage{pgfplots}
\usetikzlibrary{matrix}
\pgfplotsset{compat=newest}
\usepgfplotslibrary{groupplots}
\usetikzlibrary{pgfplots.groupplots}
\usetikzlibrary{plotmarks}
\usetikzlibrary{calc}
\usepgfplotslibrary{external}
\pgfplotsset{compat=1.3}

\usetikzlibrary{shapes.geometric}
\usepackage{fancyhdr}
\usepackage{mathtools}
\DeclarePairedDelimiter{\ceil}{\lceil}{\rceil}
\usepackage{array}
\newcolumntype{R}{>{$}r<{$}}
\newcolumntype{V}[1]{>{[\;}*{#1}{R@{\;\;}}R<{\;]}}
\usetikzlibrary{plotmarks}
\definecolor{mycolor1}{rgb}{0.20000,0.30000,0.60000}%
\definecolor{mycolor2}{rgb}{0.60000,0.20000,0.00000}%
\usepackage{pgf}
\usepackage[utf8]{inputenc}
\usetikzlibrary{arrows,automata}
\usetikzlibrary{positioning}
\usepackage{algorithm}
\usepackage[noend]{algpseudocode}

\tikzset{
    state/.style={
           rectangle,
           rounded corners,
           draw=black, very thick,
           minimum height=2em,
           inner sep=2pt,
           text centered,
           },
}

\providecommand{\norm}[1]{\left\lVert #1 \right\rVert}
\graphicspath{{./FIGS/}}

\newcommand{\ra}[1]{\renewcommand{\arraystretch}{1}}

\begin{document}

\title{\Large Projection techniques to update the truncated SVD of evolving matrices}
\author{Vassilis Kalantzis\thanks{IBM Research, Thomas J. Watson Research Center, Yorktown Heights, NY 10598, US. Email: vkal@ibm.com}
\and Georgios Kollias\thanks{IBM Research, Thomas J. Watson Research Center, Yorktown Heights, NY 10598, US. Email: gkollias@us.ibm.com}
\and Shashanka Ubaru\thanks{IBM Research, Thomas J. Watson Research Center, Yorktown Heights, NY 10598, US. Email: shashanka.ubaru@us.ibm.com}
\and Athanasios N. Nikolakopoulos\thanks{University of Minnesota, Minneapolis, MN 55455, US. Email: anikolak@umn.edu}
\and Lior Horesh\thanks{IBM Research, Thomas J. Watson Research Center, Yorktown Heights, NY 10598, US. Email: lhoresh@us.ibm.com}
\and Kenneth L. Clarkson\thanks{IBM Research, Almaden Research Center, San Jose, CA 95120, US. Email: klclarks@us.ibm.com}}

\date{}

\maketitle


\abstract{This paper considers the problem of updating the rank-$k$ 
truncated Singular Value Decomposition (SVD) of matrices subject to 
the addition of new rows and/or columns over time. Such matrix 
problems represent an important computational kernel in applications 
such as Latent Semantic Indexing and Recommender Systems. Nonetheless, 
the proposed framework is purely algebraic and 
targets general updating problems. The algorithm presented in this paper 
undertakes a projection viewpoint and focuses on building a pair of 
subspaces which approximate the linear span of the sought singular vectors 
of the updated matrix. We discuss and analyze two different choices to 
form the projection subspaces. Results on matrices from real applications 
suggest that the proposed algorithm can lead to higher accuracy, 
especially for the singular triplets associated with the largest modulus 
singular values. Several practical details and key differences with other 
approaches are also discussed.}


\section{Introduction}

This paper considers the update of the truncated SVD of a sparse 
matrix subject to additions of new rows and/or columns. More 
specifically, let $B \in \mathbb{C}^{m\times n}$ be a matrix for 
which its rank-$k$ (truncated) SVD $B_k$ is available. 
Our goal is to obtain an approximate rank-$k$ SVD $A_k$ of matrix
\begin{equation*} \label{eq200}
  A =
 \begin{pmatrix}
   B            \\[0.3em]
   E            \\[0.3em]
  \end{pmatrix},\ \ {\rm or}\ \ 
  A =
 \begin{pmatrix}
   B  & E     
  \end{pmatrix},
\end{equation*}
where $E$ denotes the matrix of newly added rows or columns. This process can 
be repeated several times, where at each instance matrix $A$ becomes matrix 
$B$ at the next level. Note that a similar problem, not explored in this paper, 
is to approximate the rank-$k$ SVD of $B$ after modifying its (non-)zero entries, 
e.g., see \cite{zha1999updating}.

Matrix problems such as the ones above hold an important role in several 
real-world applications. One such example is Latent Semantic Indexing 
(LSI) in which the truncated SVD of the current term-document matrix needs 
to be updated after a few new terms/documents have been added to the 
collection \cite{berry1995using,deerwester1990indexing,zha1999updating}. 
Another example is the update of latent-factor-based models of user-item rating 
matrices 
in top-N recommendation \cite{cremonesi2010performance,nikolakopoulos2019eigenrec,sarwar2002incremental}. 
Additional applications in geostatistical screening can be found in 
\cite[Chapter 6]{horesh2015reduced}.

The standard approach to compute $A_k$ is to disregard any previously 
available information and apply directly to $A$ an off-the-shelf, 
high-performance, SVD solver 
\cite{baglama2005augmented,hernandez2005slepc,wu2015preconditioned,halko2011finding,ubaru2019}. 
This standard approach might be feasible when the original matrix is 
updated only once or twice, however becomes increasingly impractical 
as multiple row/column updates take place over time. Therefore, it 
becomes crucial to develop algorithms which return a reasonable 
approximation of $A_k$ while taking advantage of $B_k$. Such schemes 
have already been considered extensively for the case of full SVD \cite{brand2003fast,Gu94astable,moonen1992singular} and rank-$k$ SVD 
\cite{berry1995using,sarwar2002incremental,vecharynski2014fast,zha1999updating}. 
Nonetheless, for general-purpose matrices it is rather unclear how 
to enhance their accuracy.

\subsection{Contributions.} 

\begin{enumerate}\itemsep 0pt 
    \item We propose and analyze a projection scheme to update the rank-$k$ 
    SVD of evolving matrices. Our scheme uses a right singular projection 
    subspace equal to $\mathbb{C}^n$, and only determines the left singular 
    projection subspace.
    \item We propose and analyze two different options to set the left 
    singular projection subspace. A complexity analysis is also presented.
    \item We present experiments performed on matrices stemming from 
    applications in LSI and recommender systems. These experiments demonstrate 
    the numerical behavior of the proposed scheme and showcase the various 
    tradeoffs in accuracy versus complexity.
\end{enumerate}

\section{Background and notation}

The (full) SVD of matrix $B$ is denoted as $B=U\Sigma V^H$ where 
$U\in \mathbb{C}^{m\times m}$ and $V\in \mathbb{C}^{n\times n}$ are unitary 
matrices whose $j$'th column is equal to the left singular vector $u^{(j)}$ 
and right singular vector $v^{(j)}$, respectively. The matrix $\Sigma \in \mathbb{R}^{m\times n}$ 
has non-zero entries only along its main diagonal, and these 
entries are equal to the singular values $\sigma_1\geq\cdots\geq
\sigma_{\mathtt{min}(m,n)}$. Moreover, we define the matrices 
$U_{j} = \left[u^{(1)},\ldots,u^{(j)}\right]$, 
$V_{j} = \left[v^{(1)},\ldots,v^{(j)}\right]$, 
and $\Sigma_{j} = \mathtt{diag}\left(\sigma_1,
\ldots,\sigma_j\right)$. The rank-$k$ truncated SVD of matrix $B$ 
can then be written as $B_k=U_k\Sigma_kV_k^H
=\sum_{j=1}^k \sigma_j u^{(j)} \left(v^{(j)}\right)^H$.
We follow the same notation for matrix $A$ with the exception that 
a circumflex is added on top of each variable, i.e., 
$A_k=\widehat{U}_k\widehat{\Sigma}_k\widehat{V}_k^H
=\sum_{j=1}^k \widehat{\sigma}_j \widehat{u}^{(j)} \left(\widehat{v}^{(j)}\right)^H$, 
with $\widehat{U}_{j} = \left[\widehat{u}^{(1)},\ldots,\widehat{u}^{(j)}\right]$, 
$\widehat{V}_{j} = \left[\widehat{v}^{(1)},\ldots,\widehat{v}^{(j)}\right]$, 
and $\widehat{\Sigma}_{j} = \mathtt{diag}\left(\widehat{\sigma}_1,
\ldots,\widehat{\sigma}_j\right)$.

The routines $\mathtt{nr}(K)$ and $\mathtt{nnz}(K)$ return 
the number of rows of matrix and non-zero entries of matrix 
$K$, respectively. Throughout this paper $\|\cdot\|$ will stand 
for the $\ell_2$ norm when the input is a vector, and the 
spectral norm when the input is a matrix. Moreover, the term $\mathtt{range}(K)$ will 
denote the column space of matrix $K$, while $\mathtt{span}(\cdot)$ 
will denote the linear span of a set of vectors. The identity 
matrix of size $n$ will be denoted by $I_n$.

\subsection{Related work.}

The problem of updating the SVD of an evolving matrix has been 
considered extensively in the context of LSI. Consider first 
the case 
$A = \begin{pmatrix}
   B               \\
   E               \\
\end{pmatrix}$, and let $(I-V_kV_k^H)E^H=QR$ such that $Q$ 
is orthonormal and $R$ is upper trapezoidal. The scheme in \cite{zha1999updating} 
writes
\setlength\arraycolsep{1.4pt}%
{\small \begin{eqnarray*}
\begin{pmatrix}
   B               \\[0.3em]
   E               \\[0.3em]
\end{pmatrix} \approx
\begin{pmatrix}
   U_k\Sigma_kV_k^H             \\[0.3em]
   E                            \\[0.3em]
\end{pmatrix} & = &
\begin{pmatrix}
   U_k   &           \\[0.3em]
   & I_s               \\[0.3em]
\end{pmatrix}
\begin{pmatrix}
   \Sigma_k &            \\[0.3em]
   EV_k & R^H          \\[0.3em]
\end{pmatrix}
\begin{pmatrix}
   V_k   & Q    \\[0.3em]
\end{pmatrix}^H \\
&=&
\left(\begin{pmatrix}
   U_k   &           \\[0.3em]
   & I_s               \\[0.3em]
\end{pmatrix}F\right)\Theta
\left(\begin{pmatrix}
   V_k   & Q          \\[0.3em]
\end{pmatrix}G\right)^H 
\end{eqnarray*}}
where the matrix product $F\Theta G^H$ denotes the compact 
SVD of the matrix
$\begin{pmatrix}
   \Sigma_k &            \\
   EV_k & R^H            \\
\end{pmatrix}$. 

The above idea can be also applied to $A=\begin{pmatrix}B&E\end{pmatrix}$. 
Indeed, if matrices $Q$ and $R$ are now determined as $(I-U_kU_k^H)E=QR$, 
we can approximate 
{\small \begin{eqnarray*}
\begin{pmatrix}
   B  & E
\end{pmatrix}&\approx&
\begin{pmatrix}
   U_k\Sigma_kV_k^H  &E
\end{pmatrix} \\ 
& = &
\begin{pmatrix}
   U_k   & Q
\end{pmatrix}
\begin{pmatrix}
   \Sigma_k & U_k^HE \\[0.3em]
    & R              \\[0.3em]
\end{pmatrix}
\begin{pmatrix}
   V_k^H   &         \\[0.3em]
   & I_s               \\[0.3em]
\end{pmatrix}\\
&= &
\left(\begin{pmatrix}
   U_k   & Q
\end{pmatrix}
F\right)\Theta\left(
\begin{pmatrix}
   V_k   &           \\[0.3em]
   & I_s               \\[0.3em]
\end{pmatrix}G\right)^H
\end{eqnarray*}}
where the matrix product $F\Theta G$ now denotes the compact 
SVD of the matrix
$\begin{pmatrix}
\Sigma_k & U_k^HE           \\
         & R                \\
\end{pmatrix}$.

When $B_k$ coincides with the compact SVD of $B$, the above schemes 
compute the exact rank-$k$ SVD of $A$, and no access to matrix $B$ 
is required. Nonetheless, the application of the method in 
\cite{zha1999updating} can be challenging. For general updating problems,  
or problems where $A$ does not satisfy a ``low-rank plus shift'' structure \cite{zha2000matrices}, replacing $B$ by $B_k$ might not lead to a 
satisfactory approximation of $A_k$. Moreover, the memory/computational 
cost associated with the computation of the QR and SVD decompositions in 
each one of the above two scenarios might be prohibitive. The latter was 
recognized in \cite{vecharynski2014fast} where it was proposed to 
adjust the method in \cite{zha1999updating} by replacing matrices $(I-V_kV_k^H)E^H$ and $(I-U_kU_k^H)E$ with a low-rank approximation computed by applying the Golub-Kahan Lanczos bidiagonalization 
procedure \cite{golub1965calculating}. Similar ideas have been suggested in \cite{yamazaki2017sampling} and \cite{ubaru2019sampling} where the Golub-Kahan Lanczos bidiagonalization procedure was replaced by randomized SVD \cite{halko2011finding,ubaru2015low} and 
graph coarsening \cite{ubaru2019sampling}, respectively.

\section{The projection viewpoint} \label{Sec2}

The methods discussed in the previous section can be recognized 
as instances of a Rayleigh-Ritz projection procedure and can be 
summarized as follows \cite{vecharynski2014fast,yamazaki2017sampling}:
\begin{enumerate}\itemsep 0pt
    \item Compute matrices $Z$ and $W$ such that $\mathtt{range}(Z)$ 
    and $\mathtt{range}(W^H)$ approximately capture 
    $\mathtt{range}(\widehat{U}_k)$ and $\mathtt{range}(\widehat{V}_k^H)$, 
    respectively.
    \item Compute $[\Theta_k,F_k,G_k] = \mathtt{svd}(Z^HAW)$ where 
    $\Theta_k,\ F_k$, and $G_k$ denote the $k$ leading singular values
    and associated left and right singular vectors of $Z^HAW$, respectively.
    \item Approximate $A_k$ by the product $(ZF_k)\Theta_k (WG_k)^H$.
\end{enumerate}
Ideally, the matrices $Z$ and $W$ should satisfy
\begin{eqnarray*}
    \mathtt{span}\left(\widehat{u}^{(1)},\ldots,\widehat{u}^{(k)}\right) &\subseteq& \mathtt{range}(Z),\ \ \rm{and}    \\ 
    \mathtt{span}\left(\widehat{v}^{(1)},\ldots,\widehat{v}^{(k)}\right) &\subseteq& \mathtt{range}(W).
\end{eqnarray*}
Moreover, the size of matrix $Z^HAW$ should be as small as possible to 
avoid high computational costs during the computation of $[\Theta_k,F_k,
G_k] = \mathtt{svd}(Z^HAW)$. 

Table \ref{table0} summarizes a few options to set matrices $Z$ and $W$ 
for the row updating problem. 
The method in \cite{vecharynski2014fast} considers the same matrix $Z$ as 
in \cite{zha1999updating} but sets $W=[V_k,X_r]$ where $X_r$ denotes the 
$r\in \mathbb{Z}^*$ leading left singular vectors of $(I-V_kV_k^H)E^H$.
The choice of matrices $Z$ and $W$ listed under the option “Algorithm 
\ref{alg1}" is explained in the next section. Note that the first variant 
of Algorithm \ref{alg1} uses the same $Z$ as in \cite{zha1999updating} and 
\cite{vecharynski2014fast} but different $W$. This choice leads to similar 
or higher accuracy than the scheme in \cite{zha1999updating} and this is 
also achieved asymptotically faster. A detailed comparison is deferred to 
the Supplementary Material. The second variant of Algorithm \ref{alg1} is a more 
expensive but also more accurate version of the first variant.

\begin{table}
\centering
\caption{{\it Different options to set the projection matrices $Z$ and 
$W$ for the row updating problem.}} \label{table0}
\vspace{0.05in}
\begin{tabular}{ l c c }
  \toprule
  \toprule
  Method\phantom{eigenrecAA} & $Z$ & $W$  \\
  \midrule
   \myrowcolour
  \cite{berry1995using} &  & $V_k$ \\
  \cite{zha1999updating} & {\small $Z = \begin{pmatrix}
  U_k &               \\
    & I_s              \\
\end{pmatrix}$} & 
$[V_k,Q]$\\
   \myrowcolour
  \cite{vecharynski2014fast} &  & $[V_k,X_r]$\\  \midrule
    Alg. \ref{alg1} & {\small $Z = \begin{pmatrix}
  U_k &               \\
    & I_s              \\
\end{pmatrix}$} & $I_n$ \\
\myrowcolour
  Alg. \ref{alg1} & {\small $Z = \begin{pmatrix}
  U_k,X_{\lambda,r} &               \\
    & I_s              \\
\end{pmatrix}$} & $I_n$ \\
\bottomrule
\bottomrule
\end{tabular}
\vspace*{-0.3cm}
\end{table}

\subsection{The proposed algorithm.}

Consider again the SVD update of matrix 
$A = \begin{pmatrix}
   B               \\
   E               \\
\end{pmatrix}$, with $E\in \mathbb{C}^{s\times n}$. The right 
singular vectors of $A$ trivially satisfy
$\widehat{v}^{(i)} \subseteq \mathtt{range}(I_n),\ i=1,\ldots,n$. Therefore, 
we can simply set $W=I_n$ and compute the $k$ leading singular triplets $\left(\theta_i,f^{(i)},g^{(i)}\right)$ of the matrix $Z^HAW = Z^HA$. Indeed, this 
choice of $W$ is ideal 
in terms of accuracy while it also removes the need to compute an approximate factorization of matrix $(I-V_kV_k^H)E^H$. On the other hand, the number 
of columns in matrix $Z^HAW$ is now equal to $n$ instead of $k + s$ in \cite{zha1999updating} 
and $k+l,\ l\ll s$, in \cite{vecharynski2014fast,yamazaki2017sampling}. 
This difference can be important when the full SVD of $Z^HAW$ is computed 
as in \cite{vecharynski2014fast,yamazaki2017sampling,zha1999updating}.

Our approach is to compute the singular values of $Z^HA$ in a matrix-free 
fashion while also skipping the computation of the right singular vectors 
$G_k$. 
Indeed, the matrix $G_k$ is only needed to approximate the $k$ leading 
singular vectors $\widehat{V}_k$ of $A$. Assuming that an approximation 
$\overline{U}_k$ and $\overline{\Sigma}_k$ of the matrices 
$\widehat{U}_k$ and $\widehat{\Sigma}_k$ is available, $\widehat{V}_k$ can 
be approximated as $\overline{V}_k = A^H\overline{U}_k \overline{\Sigma}_k^{-1}$.
\begin{algorithm} 
\caption{RR-SVD (“$AA^H$" version). \label{alg1}}
\begin{algorithmic}[1]
\State {\bf Input:} $B,U_k,\Sigma_k,V_k,E,Z$
\State {\bf Output:} $\overline{U}_k\approx \widehat{U}_k,\overline{\Sigma}_k\approx \widehat{\Sigma}_k,\overline{V}_k \approx \widehat{V}_k$
\State Solve $[\Theta_k,F_k]=\mathtt{svd}_k(Z^HA)$
\State Set $\overline{U}_k=ZF_k$ and $\overline{\Sigma}_k=\Theta_k$ 
\State Set $\overline{V}_k = A^H\overline{U}_k \overline{\Sigma}_k^{-1}$
\end{algorithmic}
\end{algorithm}
\vspace*{-0.3cm}
The proposed method is sketched in Algorithm \ref{alg1}. In terms of 
computational cost, Steps 4 and 5 require approximately 
$2\mathtt{nnz}(Z)k$ and $(2\mathtt{nnz}(A)+n)k$ Floating Point Operations 
(FLOPs), respectively. The complexity of Step 3 will generally depend on 
the algorithm used to compute the matrices $\Theta_k$ and $F_k$. We assume 
that these are computed by applying the unrestarted Lanczos method to 
matrix $Z^HAA^HZ$ in a matrix-free fashion \cite{saad2011numerical}. 
Under the mild assumption that Lanczos performs $\delta$ iterations for 
some $\delta \in \mathbb{Z}^*$ which is greater than or equal to $k$, 
a rough estimate of the total computational cost of Step 3 is 
$4\left(\mathtt{nnz}(Z^H)+\mathtt{nnz}(A)\right)\delta +2\mathtt{nr}(Z^H)\delta^2$ 
FLOPs. The exact complexity of Lanczos will depend on the choice of matrix $Z$. 
A detailed asymptotic analysis of the complexity of Algorithm \ref{alg1} 
and comparisons with other schemes are deferred to the Supplemental.

\begin{algorithm} 
\caption{RR-SVD (“$A^HA$" version). \label{alg2}}
\begin{algorithmic}[1]
\State {\bf Input:} $B,U_k,\Sigma_k,V_k,E,Z$
\State {\bf Output:} $\overline{U}_k\approx \widehat{U}_k,\overline{\Sigma}_k\approx \widehat{\Sigma}_k,\overline{V}_k \approx \widehat{V}_k$
\State Solve $[\Theta_k,G_k]=\mathtt{svd}_k(Z^HA^H)$ 
\State Set $\overline{V}_k=ZG_k$ and $\overline{\Sigma}_k=\Theta_k$ 
\State Set $\overline{U}_k = A\overline{V}_k \overline{\Sigma}_k^{-1}$
\end{algorithmic}
\end{algorithm}
Algorithm \ref {alg1} can be adapted to approximate $A_k$ for 
matrices of the form 
$A=\begin{pmatrix}
   B  & E    
\end{pmatrix}$. The complete procedure is summarized in Algorithm 
\ref{alg2}. Note that by combining Algorithms \ref{alg1} and 
\ref{alg2} we can approximate the $k$ leading singular triplets of 
matrices in which we add both new rows and columns.

Throughout the remainder of this paper we focus in updating the 
rank-$k$ SVD of matrix 
$A=\begin{pmatrix}
   B            \\
   E            \\
\end{pmatrix}$ by Algorithm \ref{alg1}. The discussion extends 
trivially to updates of matrix $A = \begin{pmatrix}
 B & E 
\end{pmatrix}$ by Algorithm \ref{alg2}.

\section[Building the projection matrix Z]{Building the projection matrix $Z$}

The accuracy of Step 5 in Algorithm \ref{alg1} depends on the 
accuracy of the approximate leading singular values and 
associated left singular vectors from Step 3. In turn, these 
quantities depend on how well $\mathtt{range}(Z)$ captures the 
singular vectors $\widehat{u}^{(1)},\ldots,\widehat{u}^{(k)}$ \cite{jia2001analysis,nakatsukasa2017accuracy}. Therefore, our 
focus lies in forming $Z$ such that the distance between the 
subspace $\mathtt{range}(Z)$ and the left singular vectors $\widehat{u}^{(1)},\ldots,\widehat{u}^{(k)}$ is as small as 
possible.

\subsection[Exploiting the left singular vectors of B]{Exploiting the left singular vectors of $B$.}
\label{choiceZ1}

The following proposition presents a closed-form expression of the 
$i$'th left singular vector of matrix $A=\begin{pmatrix}
   B            \\
   E            \\
\end{pmatrix}$.
\begin{proposition}\label{pro0}
The left singular vector $\widehat{u}^{(i)}$ associated with singular 
value $\widehat{\sigma}_i$ is equal to
\begin{equation*} 
\widehat{u}^{(i)}=
 \begin{pmatrix}
   -(BB^H-\widehat{\sigma}_{i}^2I_m)^{-1}BE^H\widehat{y}^{(i)}      \\[0.3em]
   \widehat{y}^{(i)}                                            \\[0.3em]
  \end{pmatrix},
\end{equation*}
where $\widehat{y}^{(i)}$ satisfies the equation
{\small \[\left[E\left(\sum\limits_{j=1}^{n} v^{(j)} \left(v^{(j)}\right)^H \dfrac{\widehat{\sigma}_{i}^2}{\widehat{\sigma}_{i}^2-\sigma_{j}^2} \right)E^H-\widehat{\sigma}_{i}^2I_s\right]\widehat{y}^{(i)}=0,\]}
and $\sigma_j=0$ for any $j=m+1,\ldots,n$ (when $n>m$).
\end{proposition}
\begin{proof}
Deferred to the Supplementary Material.
\end{proof}
The above representation of $\widehat{u}^{(i)}$ requires the solution 
of a nonlinear eigenvalue problem to compute $\widehat{y}^{(i)}$. 
Alternatively, we can express $\widehat{u}^{(i)}$ as follows.
\begin{proposition} \label{pro1}
The left singular vector $\widehat{u}^{(i)}$ associated with singular 
value $\widehat{\sigma}_i$ is equal to
\begin{equation*}
   \widehat{u}^{(i)} =
   \begin{pmatrix}
   u^{(1)},\ldots,u^{(\mathtt{min}(m,n))} & \\[0.3em]
   & I_s                    \\[0.3em]
  \end{pmatrix}
  \begin{pmatrix}
   \chi_{1,i}   \\[0.3em]
   \vdots       \\[0.3em]
   \chi_{\mathtt{min}(m,n),i}   \\[0.3em]
   \widehat{y}^{(i)}      \\[0.3em]
  \end{pmatrix}, 
\end{equation*}
where the scalars $\chi_{j,i}$ are equal to
\[\chi_{j,i}=-\left(Ev^{(j)}\right)^H\widehat{y}^{(i)}
\dfrac{\sigma_{j}}{\sigma_{j}^2-\widehat{\sigma}_{i}^2}.\]
\end{proposition}
\begin{proof}
Deferred to the Supplementary Material.
\end{proof}
Proposition \ref{pro1} suggests that setting
$Z=\begin{pmatrix}
   u^{(1)},\ldots,u^{(\mathtt{min}(m,n))} &       \\
                                          & I_s   \\
\end{pmatrix}$ 
should lead to an exact (in the absence of round-off errors) 
computation of $\widehat{u}^{(i)}$. In practice, we only have access 
to the $k$ leading left singular vectors of $B$, $u^{(1)},
\ldots,u^{(k)}$. The following proposition suggests that the 
distance between $\widehat{u}^{(i)}$ and the range space of 
$Z=\begin{pmatrix}
   u^{(1)},\ldots,u^{(k)} &       \\
                          & I_s   \\
\end{pmatrix}$
is at worst proportional to the ratio 
$\dfrac{\sigma_{k+1}}{\sigma_{k+1}^2-\widehat{\sigma}_{i}^2}$.

\begin{proposition} \label{pro2}
Let matrix $Z$ in Algorithm \ref{alg1} be defined as 
\begin{equation*}\label{Zpro2}
  Z = 
  \begin{pmatrix}
   u^{(1)},\ldots,u^{(k)} & \\[0.3em]
   & I_s                    \\[0.3em]
  \end{pmatrix}, 
\end{equation*}
and set $\gamma = O\left(\norm{E^H\widehat{y}^{(i)}}\right)$. 

Then, for any $i=1,\ldots,k$:
\begin{equation*}
\mathtt{min}_{z\in \mathtt{range}(Z)} \|\widehat{u}^{(i)}-z\| \leq 
\left|\dfrac{\gamma\sigma_{k+1}}{\sigma_{k+1}^2-\widehat{\sigma}_{i}^2}\right|.
\end{equation*}
\end{proposition}
\begin{proof}
Deferred to the Supplementary Material.
\end{proof}
Proposition \ref{pro2} implies that left singular 
vectors associated with larger singular values of $A$ are 
likely to be approximated more accurately.

\subsubsection[The structure of matrix ZH A]{The structure of matrix $Z^HA$.}

Setting the projection matrix $Z$ as in Proposition \ref{pro2} gives 
\begin{equation*}
    Z^HA = 
    \begin{pmatrix}
       V_k \Sigma_k & E^H
    \end{pmatrix}^H.
\end{equation*}
Therefore, each Matrix-Vector (MV) product with matrix $Z^HAA^HZ$ 
requires two MV products with matrices $\Sigma_k,\ V_k$ and 
$E$, for a total cost of about $4(nk+\mathtt{nnz}(E))$ FLOPs. Moreover, 
we have $\mathtt{nr}(Z^H)=s+k$, and thus a rough estimate of the 
cost of Step 3 in Algorithm \ref{alg1} is $4(nk+\mathtt{nnz}(E))\delta 
+ 2(s+k)\delta^2$ FLOPs.

\subsection{Exploiting resolvent expansions.} \label{choiceZ2}

The choice of $Z$ presented in Section \ref{choiceZ1} can compute the exact 
$A_k$ provided that the rank of $B$ is exactly $k$. Nonetheless, when the 
rank of $B$ is larger than $k$ and the singular values $\sigma_{k+1},\ldots,\sigma_{\mathtt{min}(m,n)}$ are not small, the accuracy 
of the approximate $A_k$ returned by Algorithm \ref{alg1} might be poor. This section presents an approach to enhance the projection matrix $Z$.

Recall that the top part of $\widehat{u}^{(i)}$ is equal to
$\widehat{f}^{(i)}=-(BB^H-\widehat{\sigma}_i^2I_m)^{-1}BE^H\widehat{y}^{(i)}$. 
In practice, even if we knew the unknown quantities $\widehat{\sigma}_i^2$ and 
$\widehat{y}^{(i)}$, the application of matrix $(BB^H-\widehat{\sigma}_i^2I_m)^{-1}$ 
for each $i=1,\ldots,k$, is too costly. The idea presented in this section 
considers the approximation of $(BB^H-\widehat{\sigma}_i^2I_m)^{-1},\ i=1,\ldots,k$, 
by $(BB^H-\lambda I_m)^{-1}$ for some fixed scalar $\lambda \in \mathbb{R}$.
\begin{lemma} \label{lem1}
Let 
\begin{equation*}
B(\lambda)=(I_m-U_kU_k^H)(BB^H-\lambda I_m)^{-1}
\end{equation*}
such that $\lambda > \widehat{\sigma}_k^2$. Then, we have that for 
any $i=1,\ldots,k$:
\[B(\widehat{\sigma}_i^2)
  =B(\lambda) \sum\limits_{\rho=0}^\infty 
   \left[(\widehat{\sigma}_i^2-\lambda)B(\lambda)\right]^\rho.\]
\end{lemma}
\begin{proof}
Deferred to the Supplementary Material.
\end{proof}
Clearly, the closer $\lambda$ is to $\widehat{\sigma}_i^2$, the more 
accurate the approximation in Lemma \ref{lem1} should be. We can 
now provide an expression for $\widehat{u}^{(i)}$ similar to that in 
Proposition \ref{pro1}. 

\begin{proposition} \label{pro34}
The left singular vector $\widehat{u}^{(i)}$ associated with singular 
value $\widehat{\sigma}_i$ is equal to
   \begin{eqnarray*}
   \widehat{u}^{(i)} & = &
   \begin{pmatrix}
   u^{(1)},\ldots,u^{(k)} & \\[0.3em]
   & I_s                    \\[0.3em]
  \end{pmatrix}
  \begin{pmatrix}
   \chi_{1,i}   \\[0.3em]
   \vdots       \\[0.3em]
   \chi_{k,i}   \\[0.3em]
   \widehat{y}^{(i)}      \\[0.3em]
  \end{pmatrix} \\ & & - 
  \begin{pmatrix}
   B(\lambda) \sum\limits_{\rho=0}^\infty 
   \left[(\widehat{\sigma}_i^2-\lambda)B(\lambda)\right]^\rho BE^H\widehat{y}^{(i)}  \\[0.3em]
                       \\[0.3em]
  \end{pmatrix}.
  \end{eqnarray*}
\end{proposition}
\begin{proof}
Deferred to the Supplementary Material.
\end{proof}
Proposition \ref{pro34} suggests a way to enhance the projection 
matrix $Z$ shown in Proposition \ref{pro2}. For example, we can 
approximate $B(\lambda) \sum\limits_{\rho=0}^\infty 
\left[(\widehat{\sigma}_i^2-\lambda)B(\lambda)\right]^\rho$ by 
$B(\lambda)$, which gives the following bound for the distance 
of $\widehat{u}^{(i)}$ from $\mathtt{range}(Z)$.

\begin{proposition} \label{pro35}
Let matrix $Z$ in Algorithm \ref{alg1} be defined as 
\begin{equation*}
   Z =    
   \begin{pmatrix}
   u^{(1)},\ldots,u^{(k)} & -B(\lambda)BE^H & \\[0.3em]
   & & I_s                                    \\[0.3em]
  \end{pmatrix}
\end{equation*}
and set $\gamma = O\left(\norm{E^H\widehat{y}^{(i)}}\right)$. 

Then, for any $\lambda\geq \hat{\sigma}_1^2$ and $i=1,\ldots,k$:
\begin{equation*}
\mathtt{min}_{z\in \mathtt{range}(Z)} \|\widehat{u}^{(i)}-z\| 
\leq 
\left|\dfrac{\gamma\sigma_{k+1}(\widehat{\sigma}_i^2-\lambda)}
{(\sigma_{k+1}^2-\widehat{\sigma}_{i}^2)\left(\sigma_{k+1}^2-\lambda\right)}\right|.
\end{equation*}
\end{proposition}
\begin{proof}
Deferred to the Supplementary Material.
\end{proof}
Compared to the bound shown in Proposition \ref{pro2}, 
the bound in Proposition \ref{pro35} is multiplied by  $\dfrac{\widehat{\sigma}_i^2-\lambda}{\sigma_{k+1}^2-\lambda}$. 
In practice, due to cost considerations, we choose a single 
value of $\lambda$ that is more likely to satisfy the above 
consideration, e.g., $\lambda \geq \widehat{\sigma}_1$.

\subsubsection[Computing the matrix B(lambda)BE H]{Computing the matrix $B(\lambda)BE^H$.}

The construction of matrix $Z$ shown in Lemma \ref{pro35} requires 
the computation of the matrix $-B(\lambda)BE^H$. 
The latter is equal to the matrix $X$ that satisfies the 
equation
\begin{equation}\label{eqBl1}
    -(BB^H-\lambda I_m)X  = (I_m-U_kU_k^H)BE^H.
\end{equation}
The eigenvalues of the matrix $-(BB^H-\lambda I_m)$ are equal to 
$\{\lambda-\widehat{\sigma}_i^2\}_{i=1,\ldots,m}$, and for any 
$\lambda > \widehat{\sigma}_1^2$, the matrix $-(BB^H-\lambda I_m)$ is 
positive definite. It is thus possible to compute $X$ by repeated 
applications of the Conjugate Gradient method. 

\begin{proposition} \label{pro46}
Let $K=-(BB^H-\lambda I_m)$ and $\|e_j\|_K$ denote the $K$-norm of 
the error after $j$ iterations of the Conjugate Gradient method 
applied to the linear system
$-(BB^H-\lambda I_m)x = b$, where $b\in \mathtt{range}
\left((I_m-U_kU_k^H)BE^H\right)$. Then,
\begin{equation*}
    \|e_j\|_K \leq 2 
    \left(\dfrac{\sqrt{\kappa}-1}{\sqrt{\kappa}+1}\right)^j \|e_0\|_K, 
\end{equation*}
where 
$\kappa = 
\dfrac{\sigma_{\mathtt{min}(m,n)}^2-\lambda}{\sigma_{k+1}^2-\lambda}$ 
and $\lambda > \widehat{\sigma}_1^2$.
\end{proposition}
\begin{proof}
Since $b\in \mathtt{range}((I_m-U_kU_k^H)BE^H)$, the vector $x$ 
satisfies the equation 
\begin{equation}\label{eqBl2}
    -(I_m-U_kU_k^H)(BB^H-\lambda I_m)(I_m-U_kU_k^H)x = b.
\end{equation}
The proof can then be found in \cite{saad2000deflated}. 
\end{proof}

\begin{corollary}
The effective condition number satisfies the inequality $\kappa \leq \dfrac{\lambda}{\lambda-\sigma_{k+1}^2}$.
\end{corollary}
Proposition \ref{pro46} applies to each one of the $s$ right-hand sides 
in (\ref{eqBl1}). Assuming that the matrix $(I_m-U_kU_k^H)BE^H$ can be 
formed and stored, the effective condition number can be reduced even 
further. For example, solving (\ref{eqBl1}) by the block Conjugate Gradient 
method leads to an effective condition number 
$\kappa=\dfrac{\sigma_{\mathtt{min}(m,n)}^2-\lambda}{\sigma_{k+s+1}^2-\lambda}$ 
\cite{o1980block}.
Additional techniques to solve linear systems with multiple right-hand 
sides can be found in 
\cite{kalantzis2013accelerating,kalantzis2018scalable,stathopoulos2010computing}.

Finally, notice that as $\lambda$ increases, the effective condition 
number decreases. Thus from a convergence viewpoint, it is better to 
choose $\lambda\gg \widehat{\sigma}_{i}^2$. On the other hand, increasing 
$\lambda$ leads to worse bounds in Proposition \ref{pro35}. 

\subsection[Truncating the matrix B(lambda)BE H]{Truncating the matrix $B(\lambda)BE^H$.} \label{trunc}

When the number of right-hand sides in (\ref{eqBl1}), i.e., 
number of rows in matrix $E$, is too large, an alternative is 
to consider $-B(\lambda)BE^H\approx X_{\lambda,r}S_{\lambda,r}
Y_{\lambda,r}^H$, where $X_{\lambda,r}S_{\lambda,r}Y_{\lambda,r}^H$ 
denotes the rank-$r$ truncated SVD of matrix $-B(\lambda)BE^H$. 
We can then replace $-B(\lambda)BE^H$ by $X_{\lambda,r}$, since $\mathtt{range}\left(X_{\lambda,r}S_{\lambda,r}Y_{\lambda,r}^H\right)
\subseteq\mathtt{range}\left(X_{\lambda,r}\right)$. 

The matrix $X_{\lambda,r}$ can be approximated in a matrix-free fashion 
by applying a few iterations of Lanczos bidiagonalization to matrix 
$B(\lambda)BE^H$. Each iteration requires two applications of Conjugate 
Gradient to solve linear systems of the same form as in (\ref{eqBl2}). 
A second approach is to apply randomized SVD as described in 
\cite{halko2011finding,clarkson2009numerical}. In practice, this amounts 
to computing the SVD of the matrix $B(\lambda)BE^HEB^HB(\lambda)R$ where $R$ is a real matrix with at least $r$ columns whose entries are i.i.d. Gaussian random variables of zero mean and unit variance.

\subsubsection[The structure of matrix \texorpdfstring{$Z^HA$}.]{The structure of matrix $Z^HA$.}

Setting the basis matrix $Z$ as in Proposition \ref{pro35}
leads to  
\begin{equation*}
    Z^HA = 
    \begin{pmatrix}
        V_k \Sigma_k &  B^HX_{\lambda,r} & E^H 
    \end{pmatrix}^H.
\end{equation*}
Each MV product with matrix $Z^HAA^HZ$ then requires two MV products 
with matrices $\Sigma_k,\ V_k,\ E$ and $B^HX_{\lambda,r}$, for a total 
cost of $4(n(k+r)+\mathtt{nnz}(E))$. Moreover, we have $\mathtt{nr}
(Z^H)=k+r+s$, and thus a rough estimate of the cost of Step 3 in Algorithm 
\ref{alg1} is $4(n(k+r)+\mathtt{nnz}(E))\delta + 2(s+k+r)\delta^2$ FLOPs.

\section{Evaluation}

Our experiments were conducted in a Matlab environment (version R2020a), 
using 64-bit arithmetic, on a single core of a computing system equipped 
with an Intel Haswell E5-2680v3 processor and 32 GB of system memory.

\begin{table}[ht]
\centering
\caption{\it Properties of the test matrices used throughout this section. 
\label{table1}}\vspace{0.05in}
\begin{tabular}{ l c c c c}
  \toprule
  \toprule
  Matrix & rows & columns
  & $nnz(A)$/rows & Source \\
  \midrule
  \myrowcolour
  MED & 5,735& 1,033 & 8.9 & \cite{berrydata} \\
  CRAN & 4,563 & 1,398 & 17.8 & \cite{berrydata} \\
  \myrowcolour
  CISI & 5,544 & 1,460 & 12.2 & \cite{berrydata} \\
  ML1M & 6,040 & 3,952 & 165.6 & \cite{harper2015movielens}\\
  \bottomrule
  \bottomrule
\end{tabular}
\end{table}

\begin{figure}[ht]
  \centering
 \includegraphics[width=0.89\linewidth]{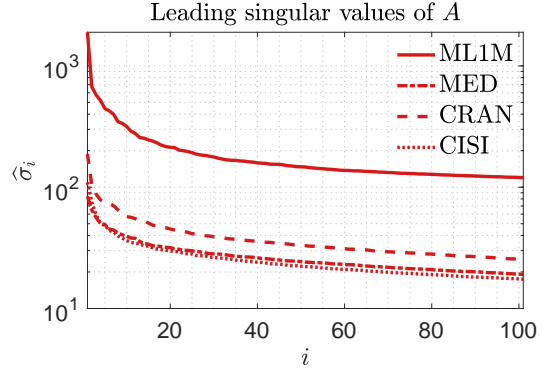}
  \caption{{\it Leading $k=100$ singular values.}}\label{fig:15}
\end{figure}

\begin{figure*}[!ht]
     \centering
     \includegraphics[width=0.24\linewidth]{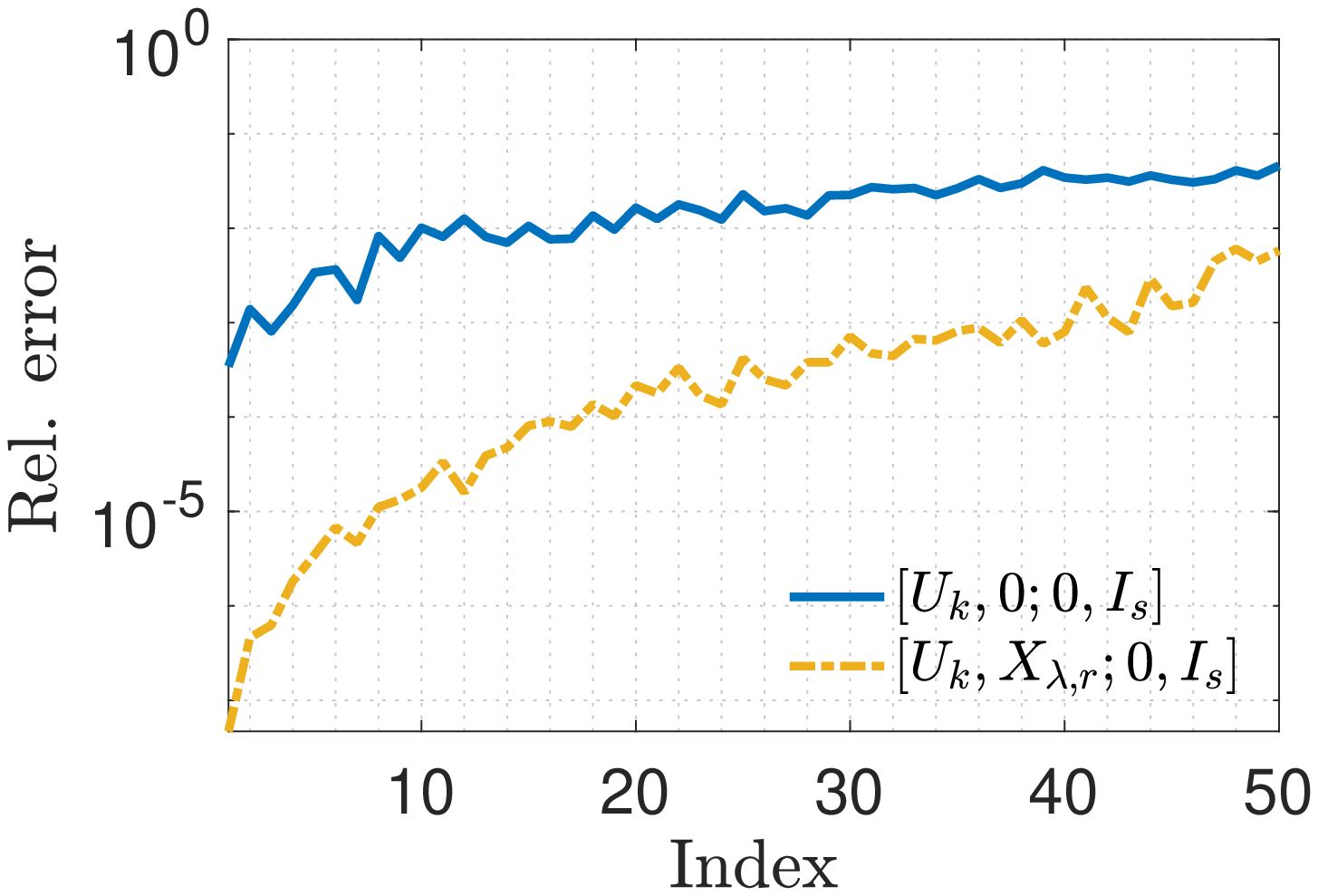}
     \includegraphics[width=0.24\linewidth]{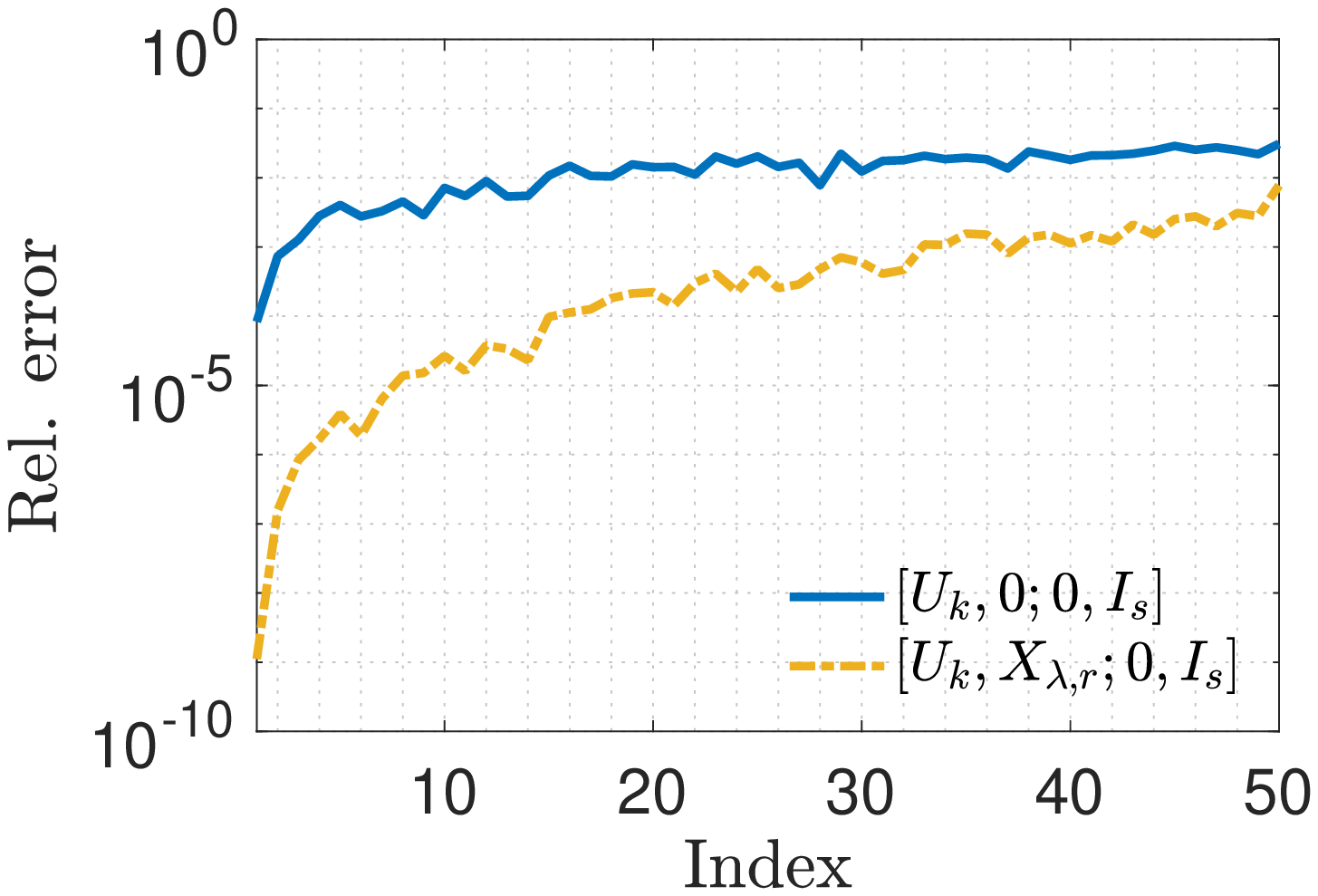}
     \includegraphics[width=0.24\linewidth]{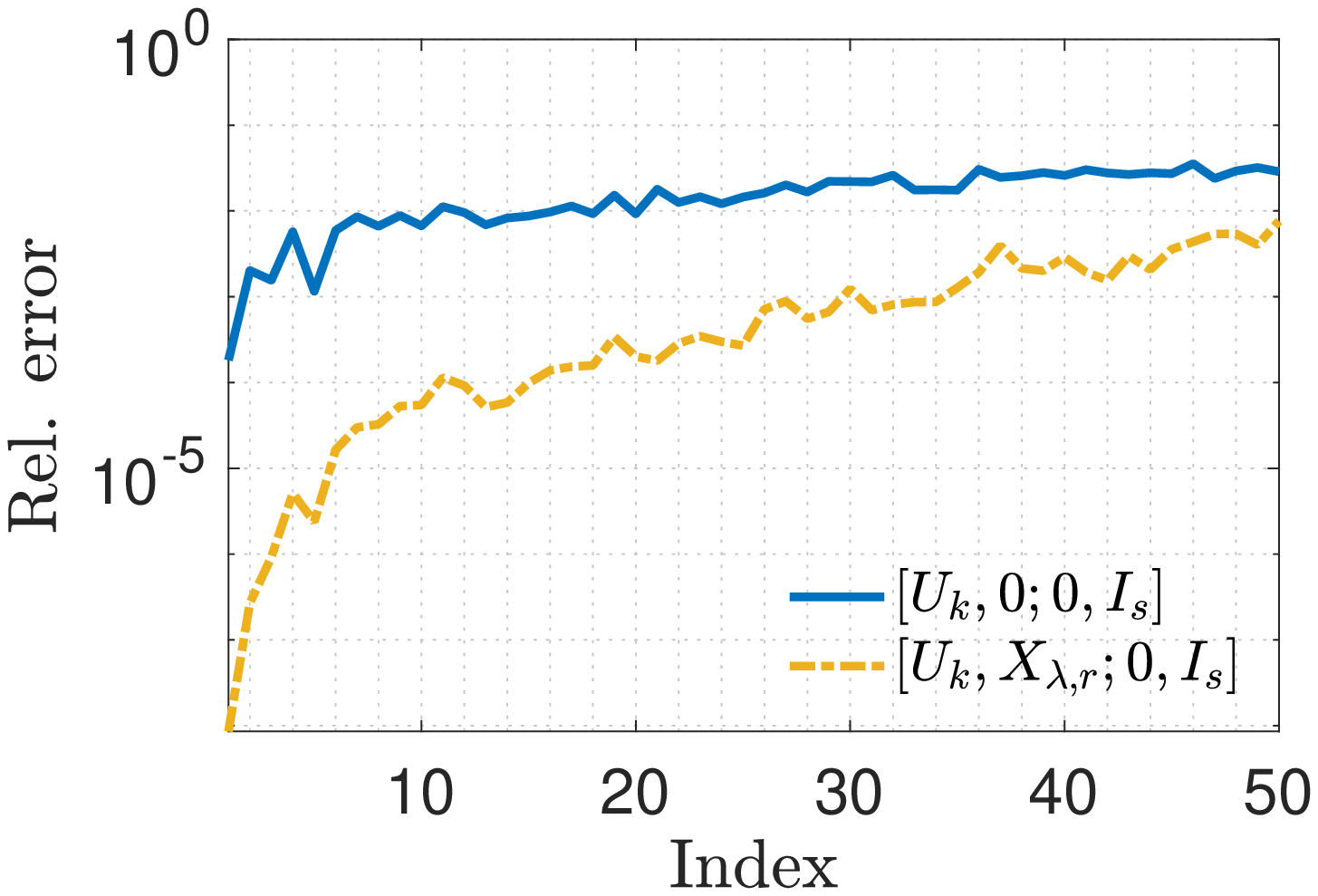}
     \includegraphics[width=0.24\linewidth]{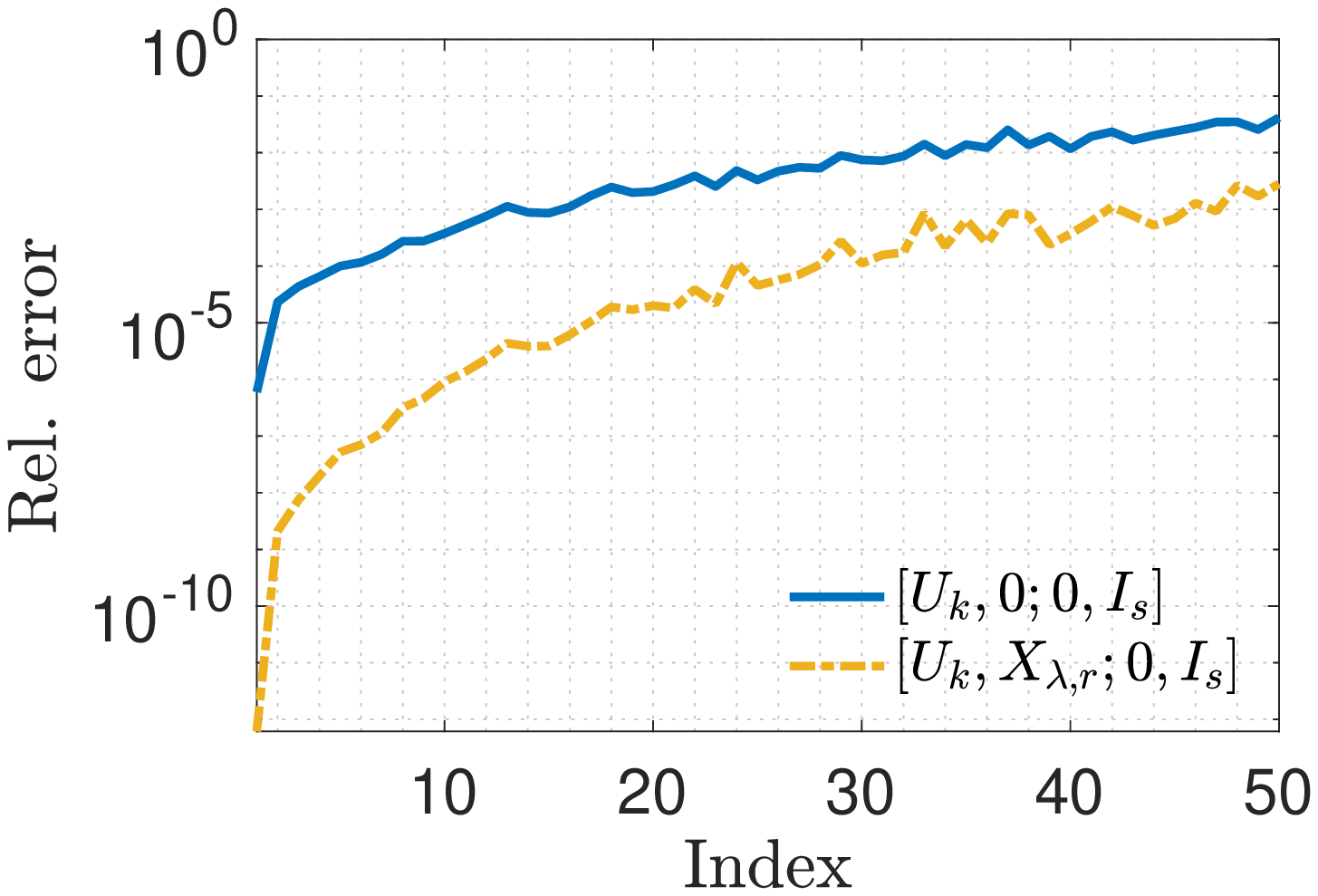}
     \includegraphics[width=0.24\linewidth]{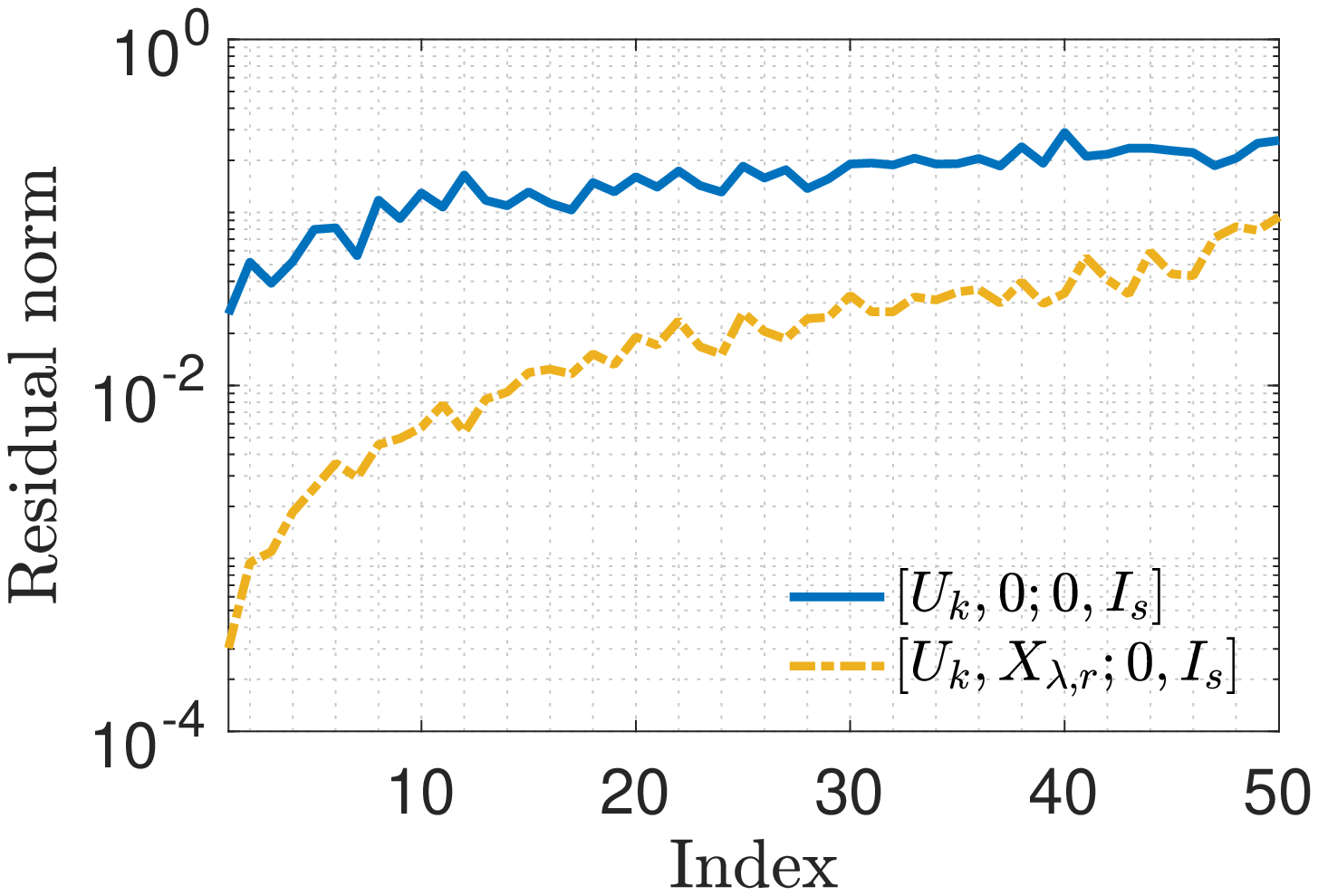}
     \includegraphics[width=0.24\linewidth]{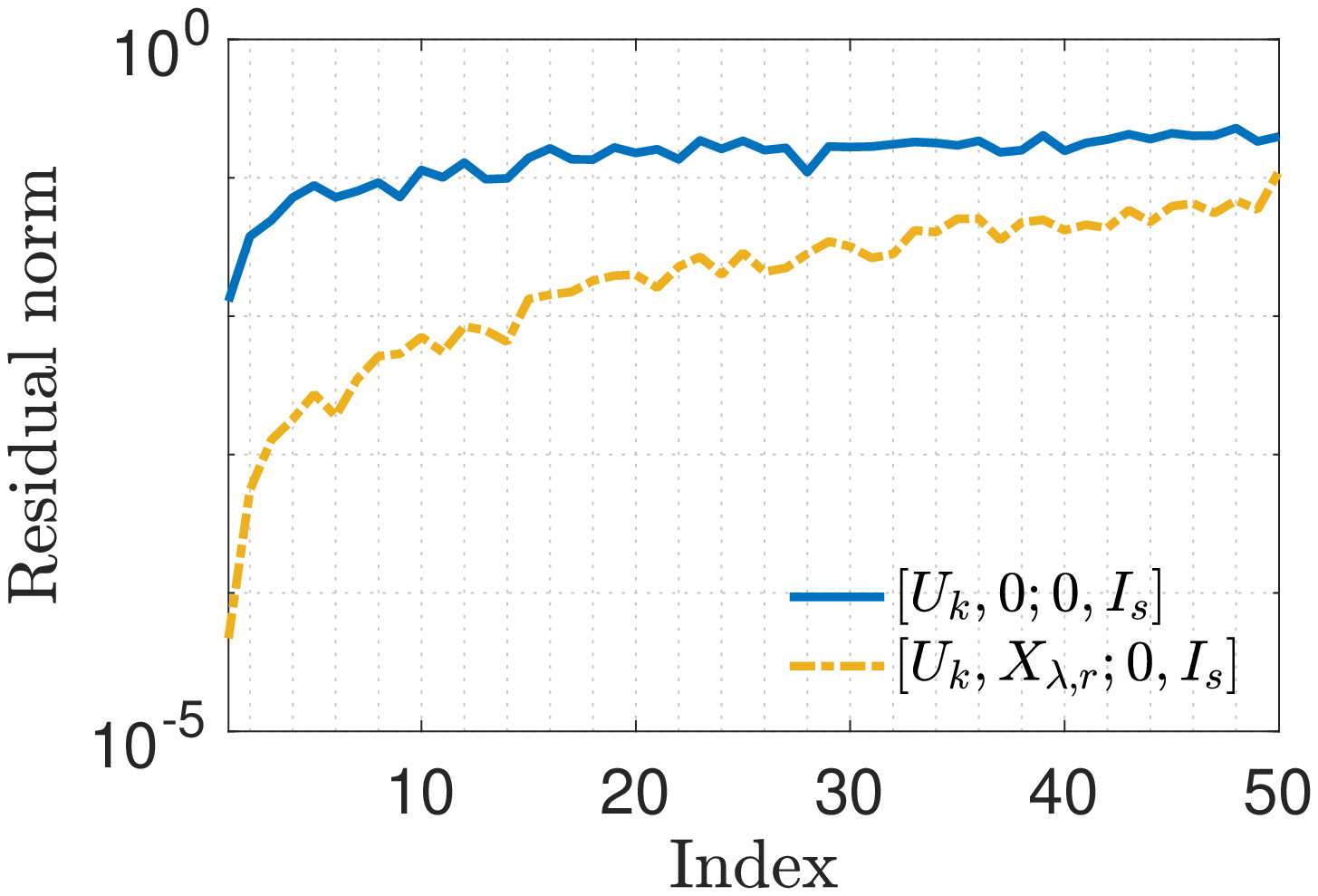}
     \includegraphics[width=0.24\linewidth]{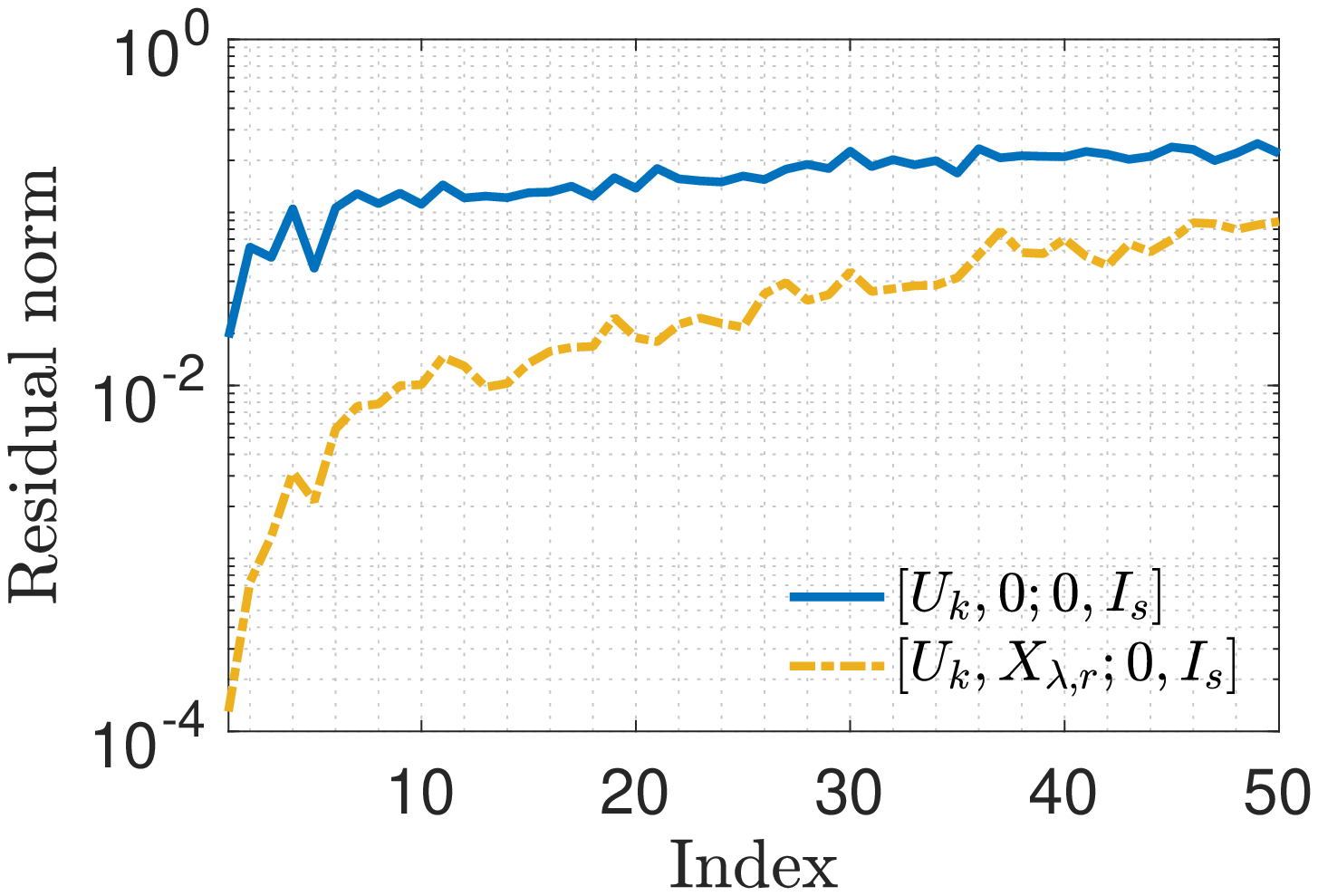}
     \includegraphics[width=0.24\linewidth]{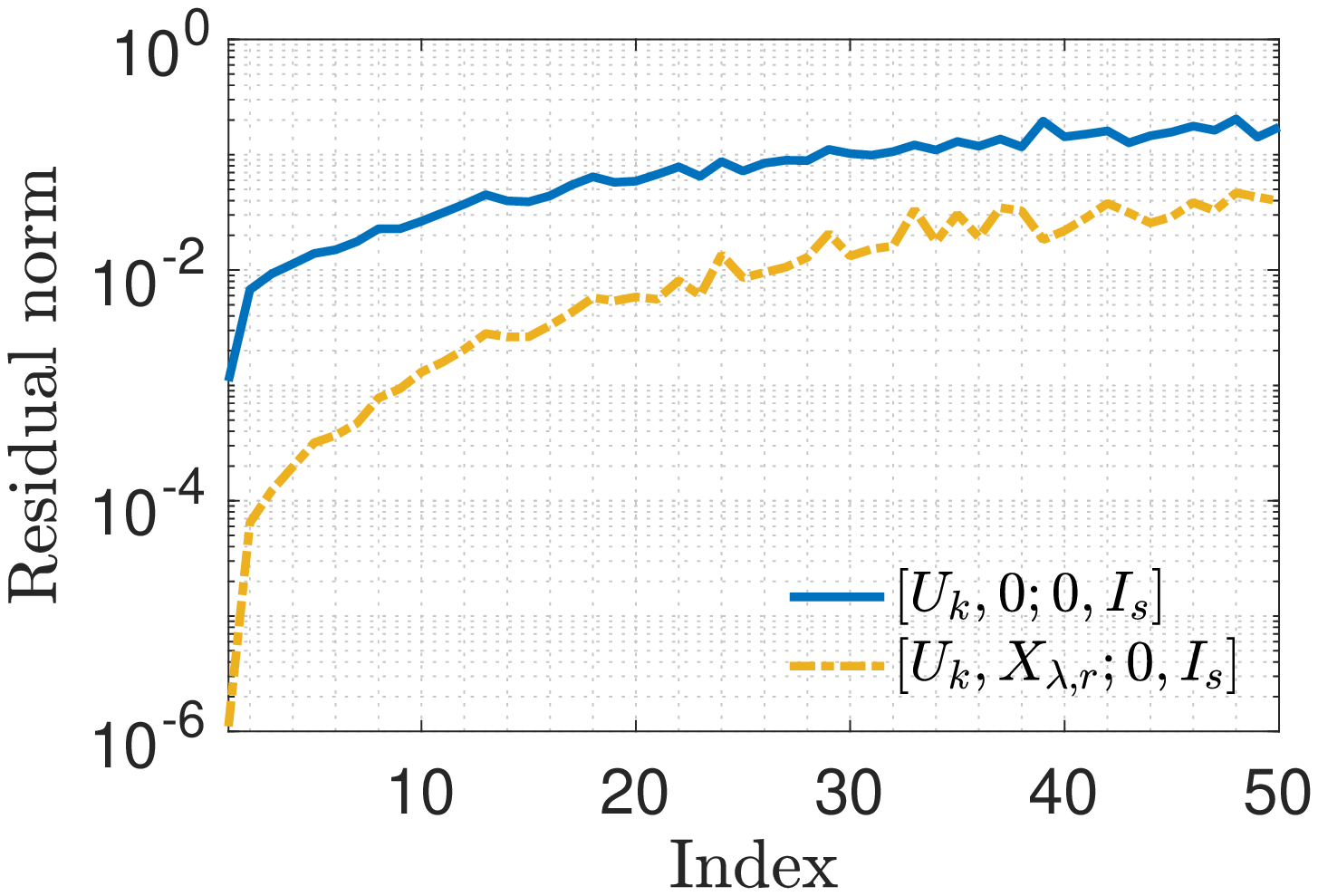}
     \caption{{\it Approximation of the leading $k=50$ singular triplets 
     for the single update case. From left to right: MED, CRAN, CISI, and 1M.}}\label{fig:16}
\end{figure*}

Table \ref{table1} lists the test matrices considered throughout 
our experiments along with their dimensions and source from which 
they were retrieved. The first three matrices come from LSI 
applications and represent term-document matrices, while the last 
matrix comes from recommender systems and represents a user-item 
rating matrix. The $k=100$ leading singular values of each matrix 
listed in Table \ref{table1} are plotted in Figure \ref{fig:15}. 

Throughout this section we focus on accuracy and will be reporting: 
a) the relative error in the approximation of the $k$ leading 
singular values of $A$, and b) the norm of the residual 
$A\widehat{v}^{(i)}-\widehat{\sigma}_i\widehat{u}^{(i)}$, scaled 
by $\widehat{\sigma}_i$. The scalar $\lambda$ is set as 
$\lambda = 1.01 \widehat{\sigma}_1^2$ where the latter singular 
value is approximated by a few iterations of Lanczos bidiagonalization.

\subsection{Single update.}

In this section we consider the approximation of the $k=50$ leading 
singular triplets of $A = \begin{pmatrix}
   B               \\
   E               \\
\end{pmatrix}$ where $B=A(1:\ceil{m/2},:)$, i.e., the size of matrices 
$B$ and $E$ is about half the size of $A$. We run Algorithm \ref{alg1} 
and set $Z$ as in Propositions \ref{pro2} and \ref{pro35}. For the 
enhanced matrix $Z$, the matrix $X_{\lambda,r}$ is computed by randomized 
SVD where $r=k$ and the number of columns in matrix $R$ is equal to $2k$ 
(recall the discussion in Section \ref{trunc}). The associated linear 
system with $2k$ right-hand sides is solved by block Conjugate Gradient. 

Figure \ref{fig:16} plots the relative error and residual norm in the 
approximation of the $k=50$ leading singular triplets of $A$. As expected, 
enhancing the projection matrix $Z$ by $X_{\lambda,r}$ leads to higher 
accuracy. This is especially true for the approximation of those singular 
triplets with corresponding singular values $\widehat{\sigma}_i\approx 
\lambda$.

In all of our experiments, the worst-case (maximum) relative error and 
residual norm was achieved in the approximation of the singular 
triplet $(\hat{\sigma}_{50},\hat{u}^{(50)},\hat{v}^{(50)})$. 
Table \ref{table2} lists the relative error and residual norm associated 
with the approximation of the singular triplet $(\hat{\sigma}_{50},\hat{u}^{(50)},\hat{v}^{(50)})$ as $r$ varies from 
ten to fifty in increments of ten. As a reference, we list the same 
quantity for the case 
{\footnotesize $Z =    
   \begin{pmatrix}
   U_k & \\
   & I_s \\
  \end{pmatrix}$}. 
As expected, enhancing the projection matrix $Z$ by $X_{\lambda,r}$ leads to higher 
accuracy, especially for higher values of $r$. 

\begin{table*}[t]
\centering
\caption{\it Relative error and residual norm associated 
with the approximation of the singular triplet $(\widehat{\sigma}_{50},\widehat{u}^{(50)},\widehat{v}^{(50)})$. \label{table2}}

\vspace{0.05in}
\begin{tabular}{l @{\hskip 0.3in} c  c c  c  c c c c c c c c}
\toprule
\toprule
\multirow{2}{*}{} & 
& \multicolumn{2}{c}{\textbf{MED}} & &
\multicolumn{2}{c}{\textbf{CRAN}} & & \multicolumn{2}{c}{\textbf{CISI}} & & \multicolumn{2}{c}{\textbf{ML1M}}\\
\cmidrule[0.4pt](lr{0.125em}){3-4}%
\cmidrule[0.4pt](lr{0.125em}){6-7}%
\cmidrule[0.4pt](lr{0.125em}){9-10}%
\cmidrule[0.4pt](lr{0.125em}){12-13}%
&  $r$ & err. & res. & & err. & res. & & err. & res. & & err. & res. \\
\midrule
\multirow{5}{*}{{\footnotesize $Z =    
   \begin{pmatrix}
   U_k & X_{\lambda,r} & \\
   & & I_s               \\
  \end{pmatrix}$}} & \tikzmarkin[hor=style mygrey]{r10}$r=10$ & 0.036 & 0.234 &  & 0.026 & 0.176 &  & 0.025 & 0.214  &  & 0.031 & 0.156\tikzmarkend{r10} \\
 & $r=20$ & 0.031 & 0.184 &  & 0.021 & 0.155  &  & 0.023 & 0.189 &  & 0.012 & 0.143 \\
  & \tikzmarkin[hor=style mygrey]{r30}$r=30$ & 0.021 & 0.114 &  & 0.017 & 0.134 &  & 0.017 & 0.161 &  & 0.008& 0.121\tikzmarkend{r30}\\
 & $r=40$ & 0.009 & 0.091 &  & 0.013 & 0.111  &  & 0.012 & 0.134 & &  0.005 & 0.112 \\
  & \tikzmarkin[hor=style mygrey]{r50}$r=50$ & 0.004 & 0.053 &  & 0.007 & 0.098 &   & 0.007 & 0.081  &  & 0.003 & 0.076\tikzmarkend{r50} \\
\midrule
{\footnotesize $Z =    
   \begin{pmatrix}
   U_k & \\
   & I_s \\
  \end{pmatrix}$}  & N/A & 0.045 & 0.269  &  & 0.045 & 0.199  &  & 0.287 & 0.250  &  & 0.041& 0.173\\
  \bottomrule
  \bottomrule
\end{tabular}
\end{table*}

\begin{figure*}[ht]
     \centering
     \includegraphics[width=0.24\linewidth]{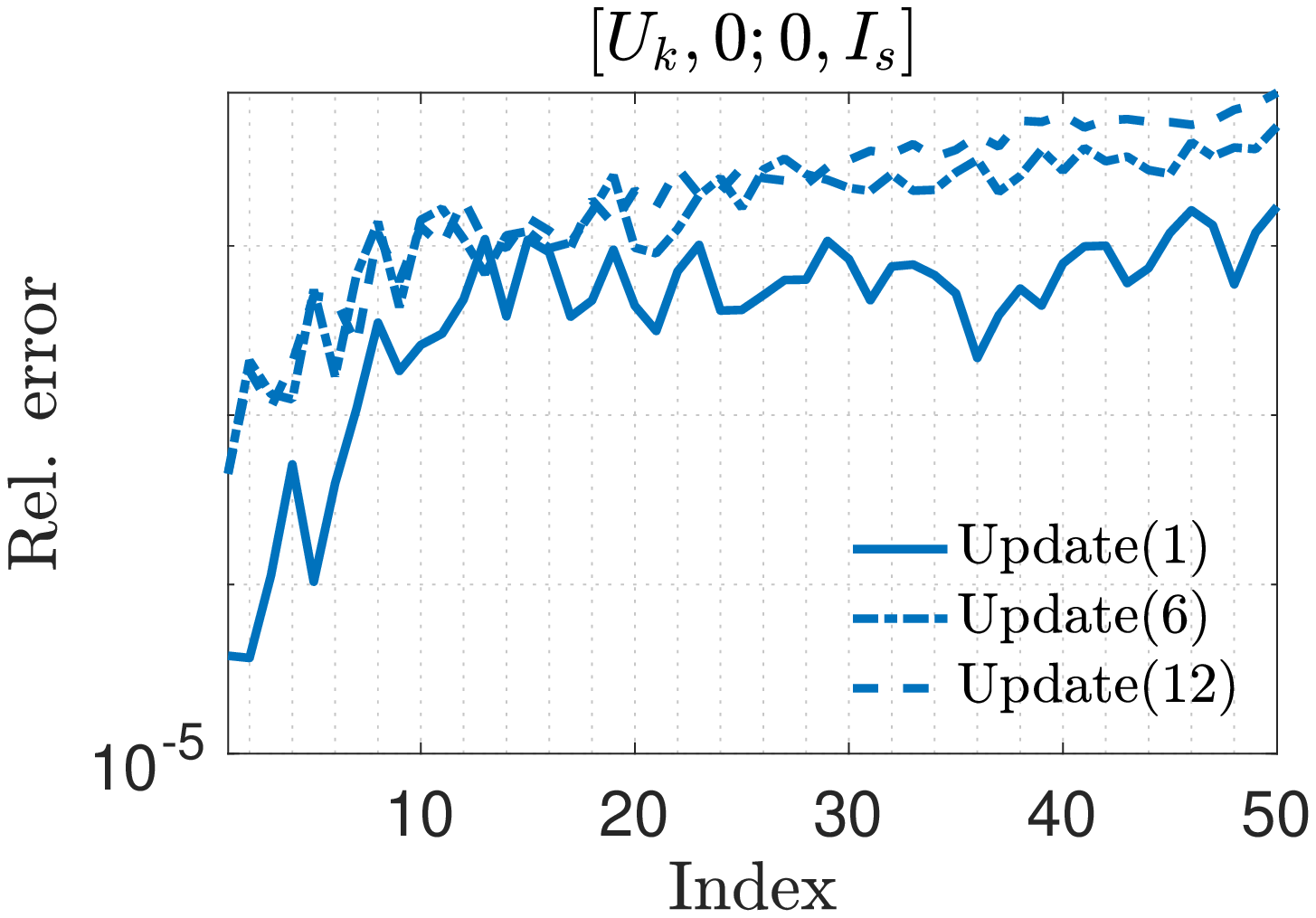}
     \includegraphics[width=0.24\linewidth]{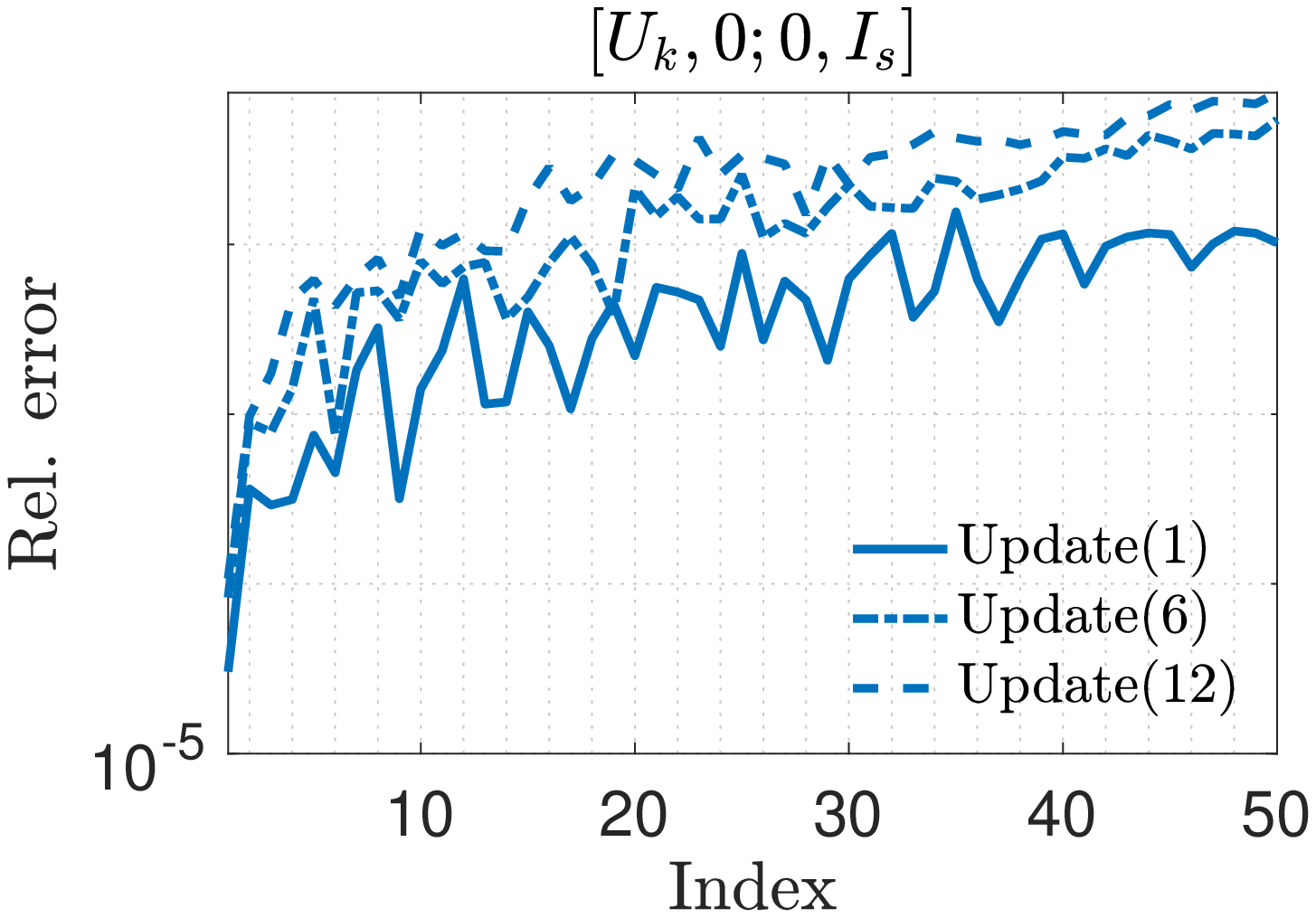}
     \includegraphics[width=0.24\linewidth]{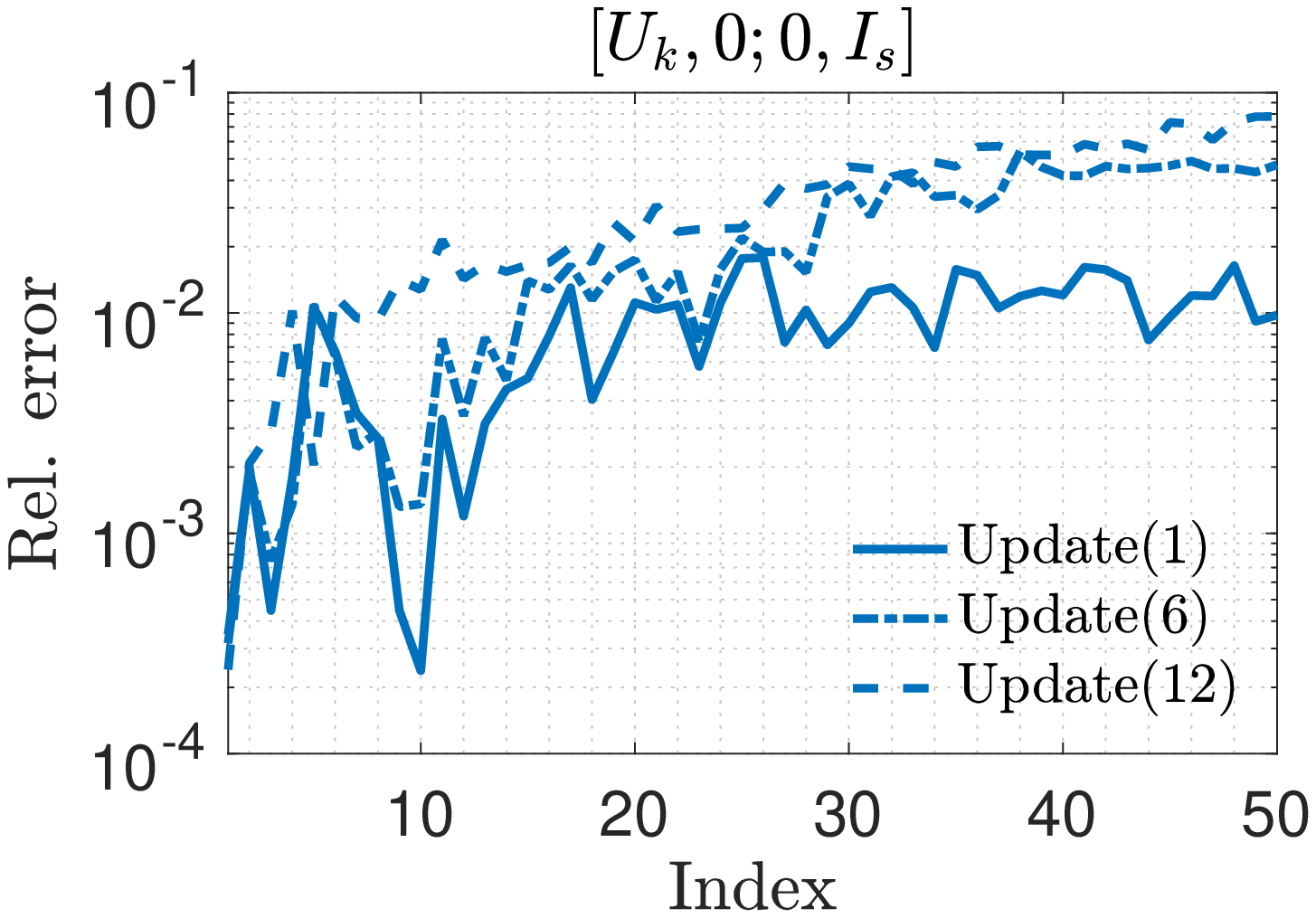}
     \includegraphics[width=0.24\linewidth]{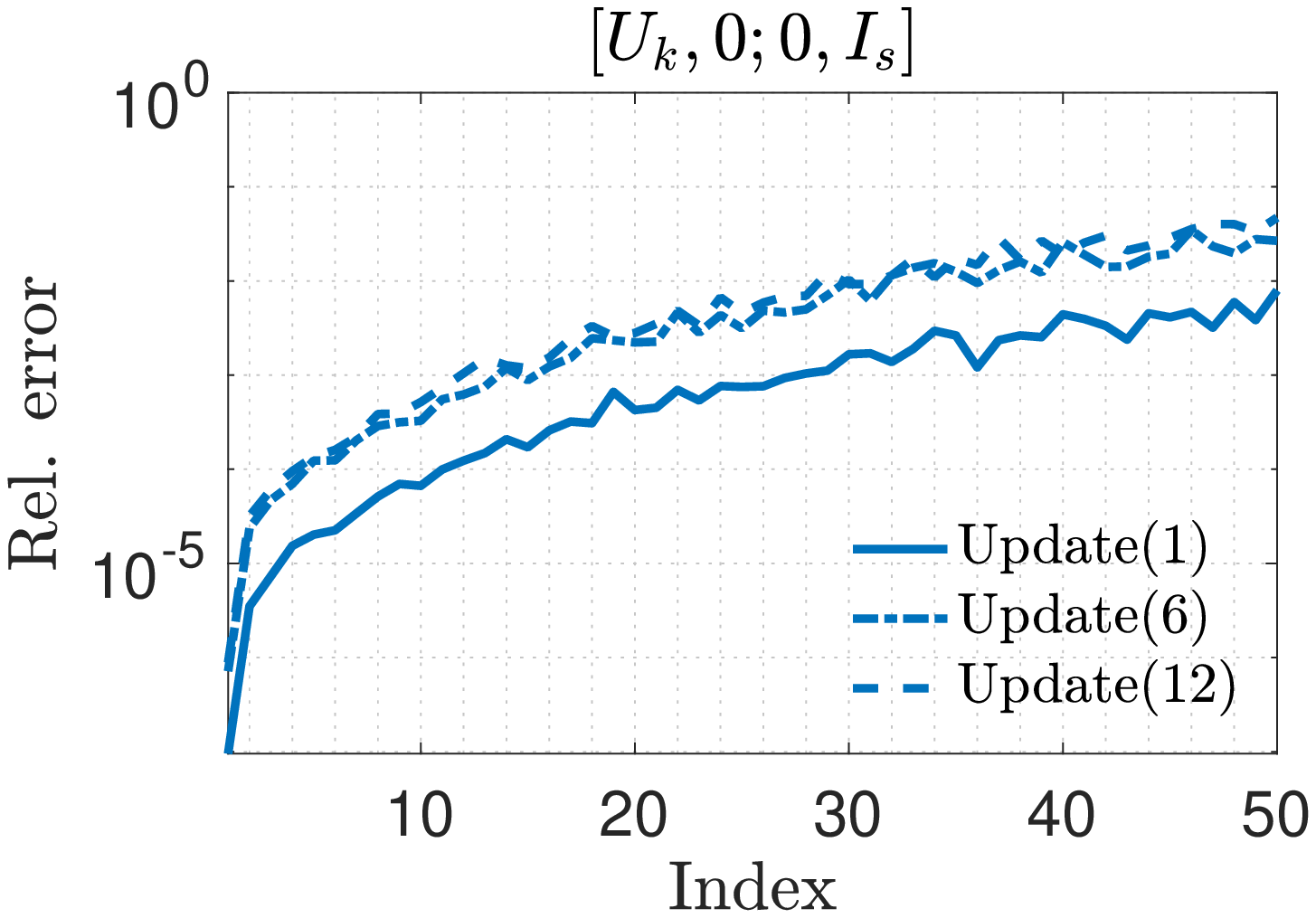}
     \includegraphics[width=0.24\linewidth]{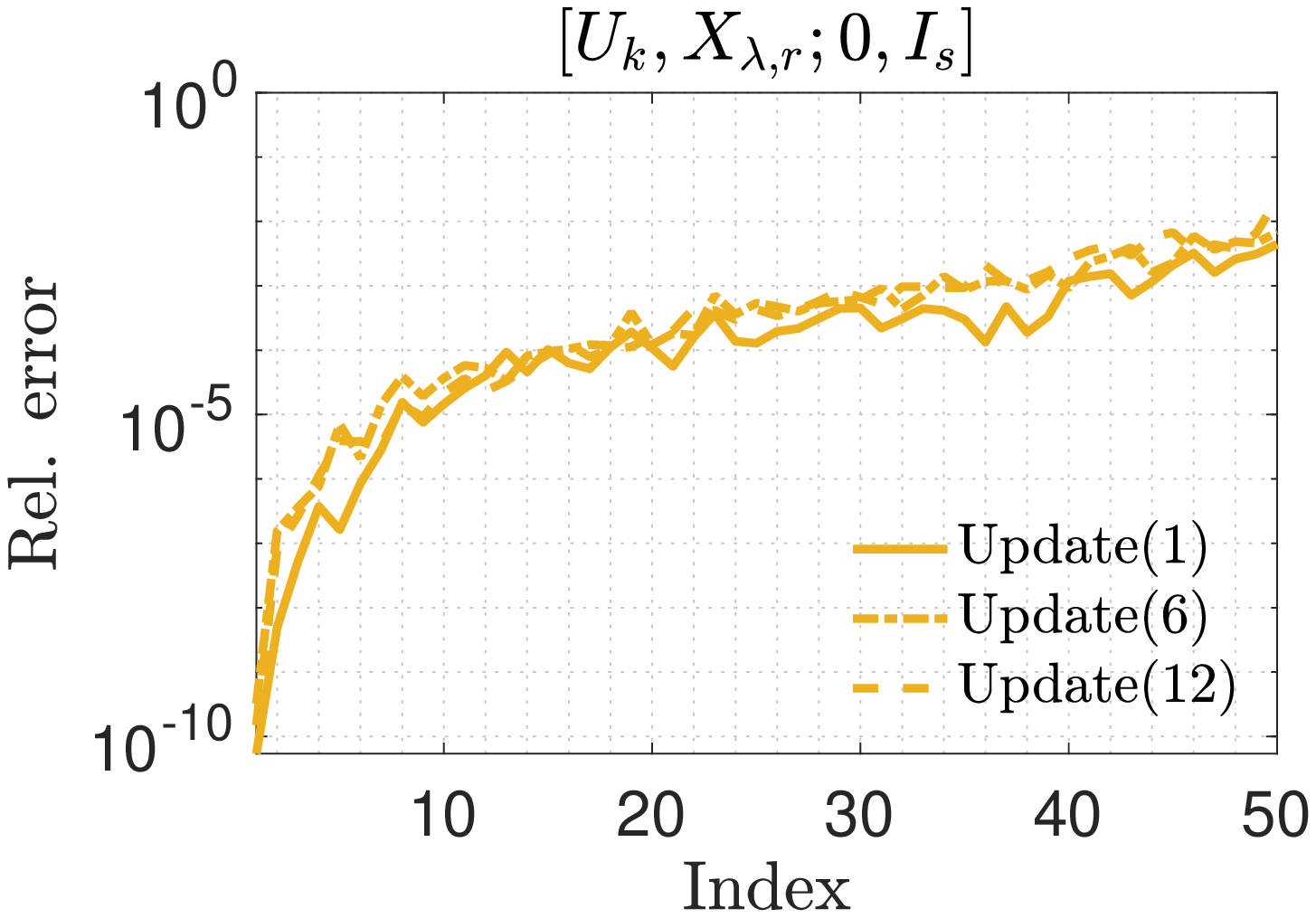}
     \includegraphics[width=0.24\linewidth]{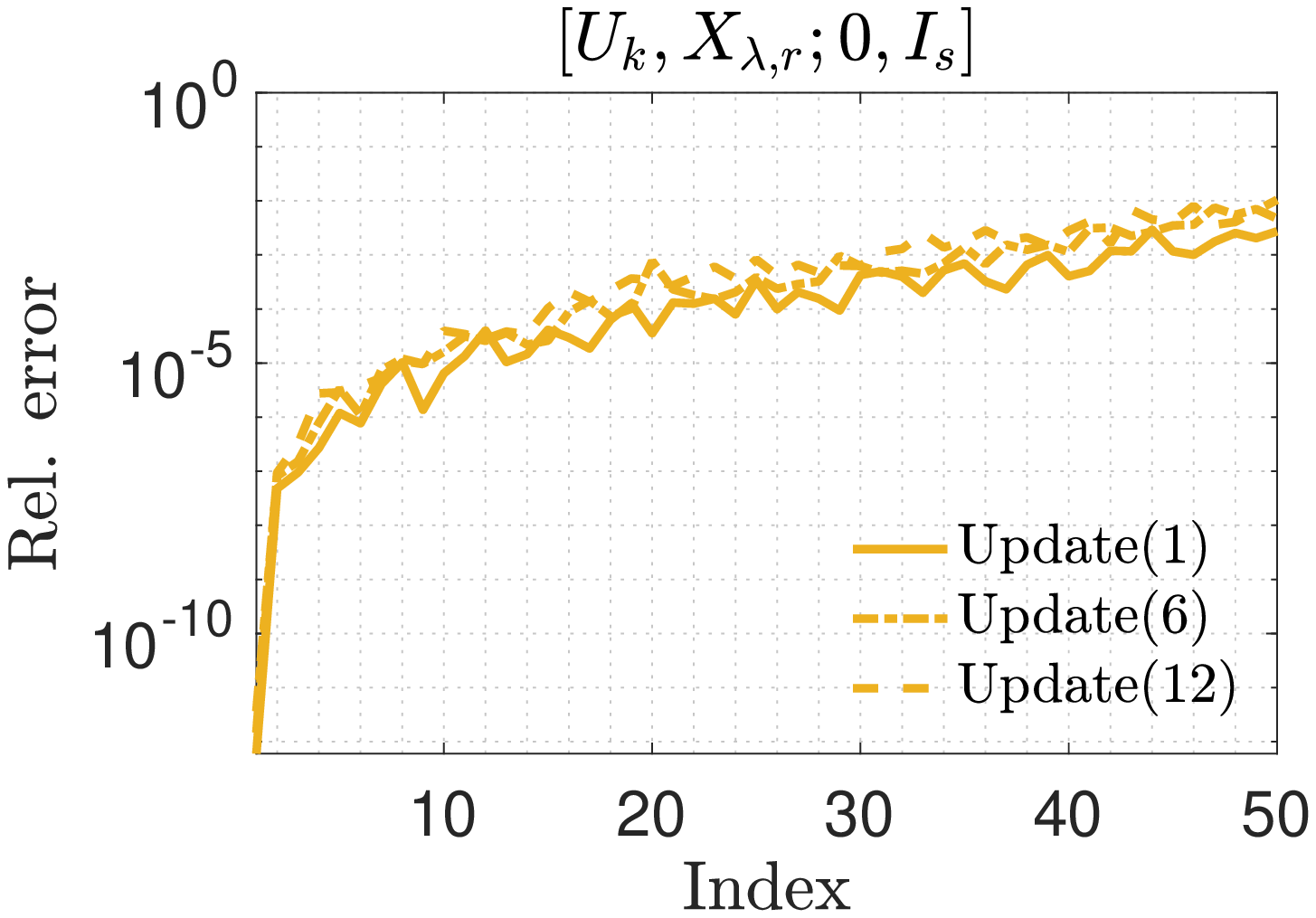}
     \includegraphics[width=0.24\linewidth]{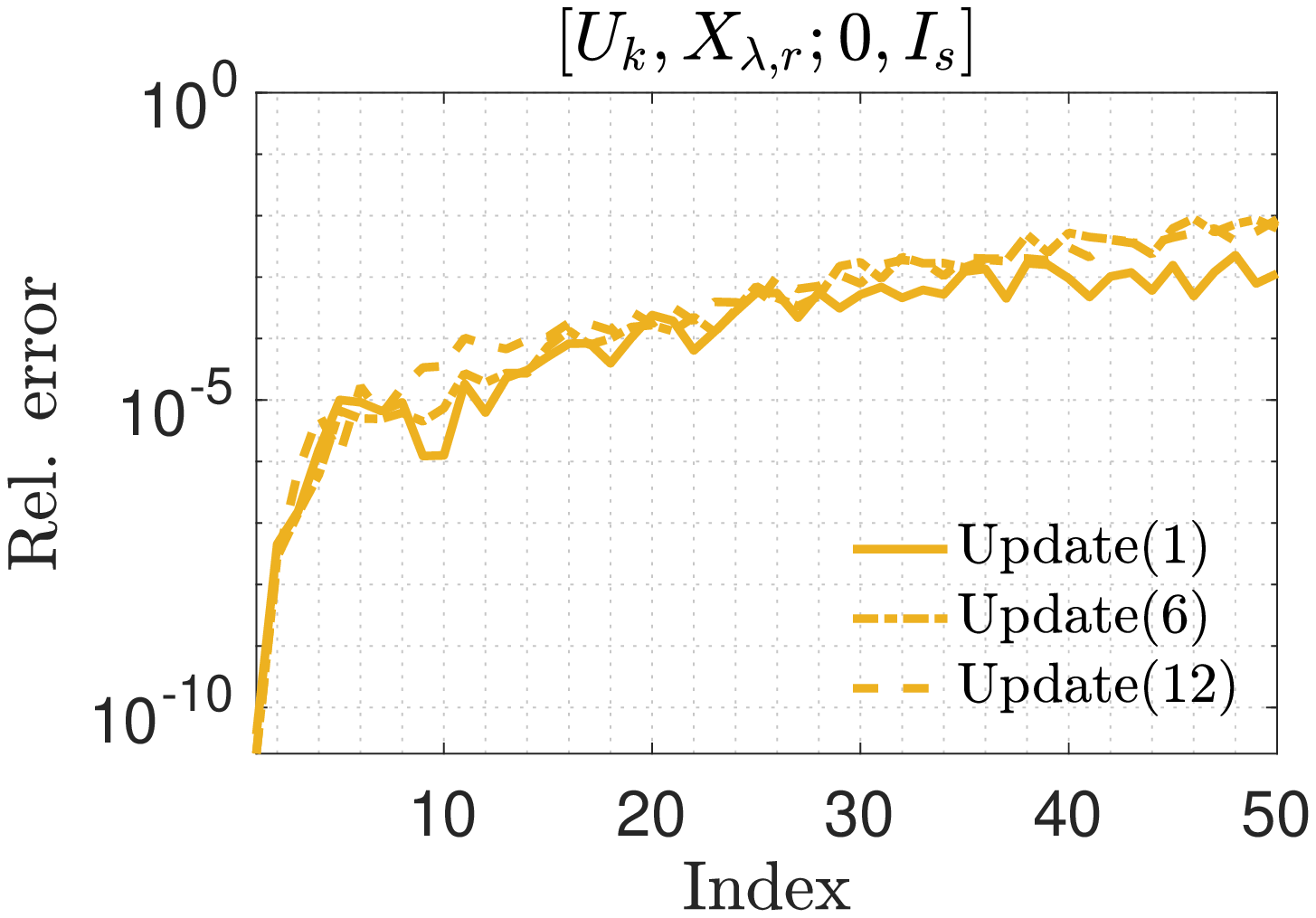}
     \includegraphics[width=0.24\linewidth]{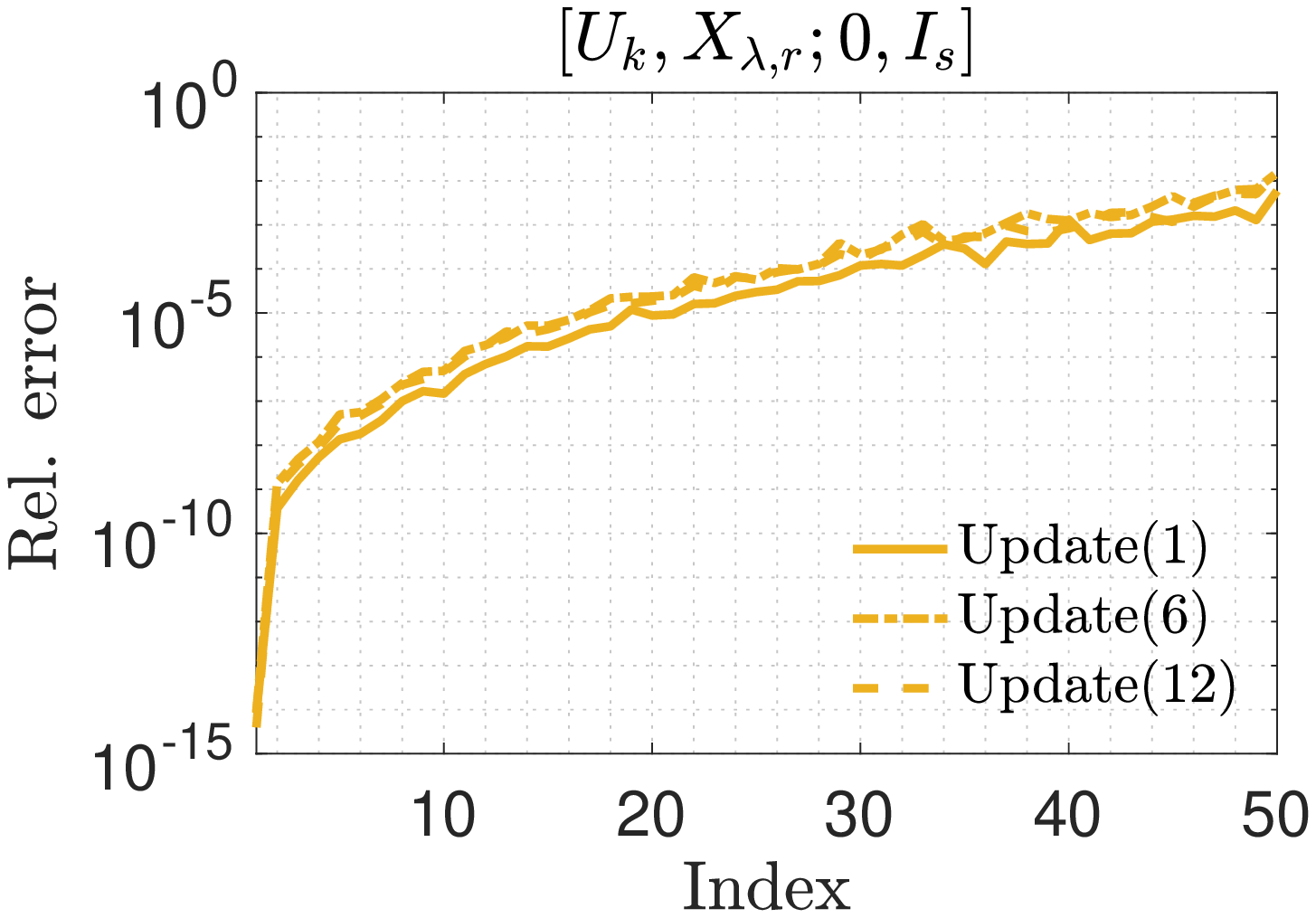}
     \caption{{\it Relative error in the approximation of the $k=50$ leading 
     singular values of $A$ for the multiple updates case. From left to right: 
     MED, CRAN, CISI, and 1M.}}\label{fig:17}
\end{figure*}

\begin{figure*}[ht]
     \centering
     \includegraphics[width=0.24\linewidth]{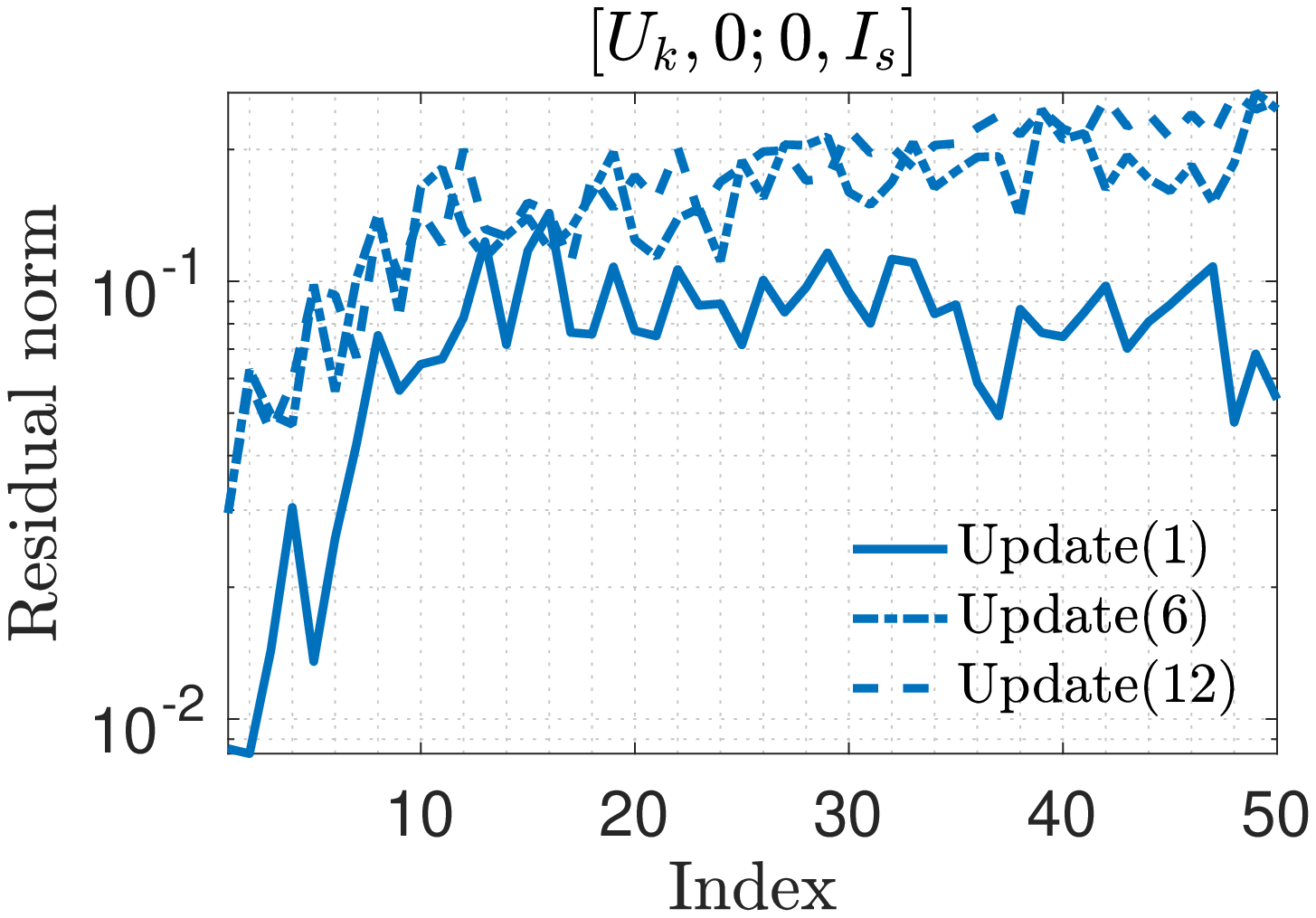}
     \includegraphics[width=0.24\linewidth]{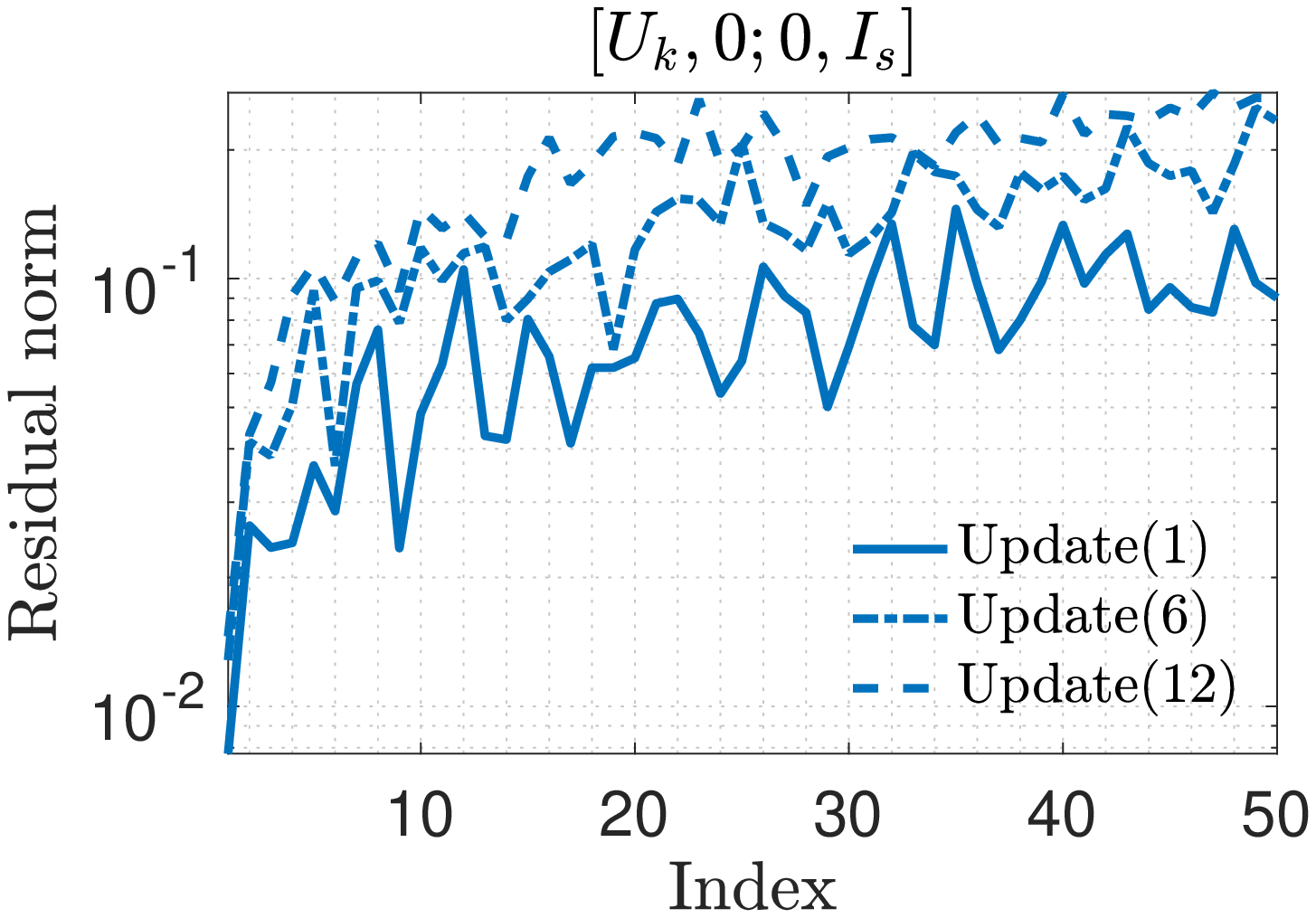}
     \includegraphics[width=0.24\linewidth]{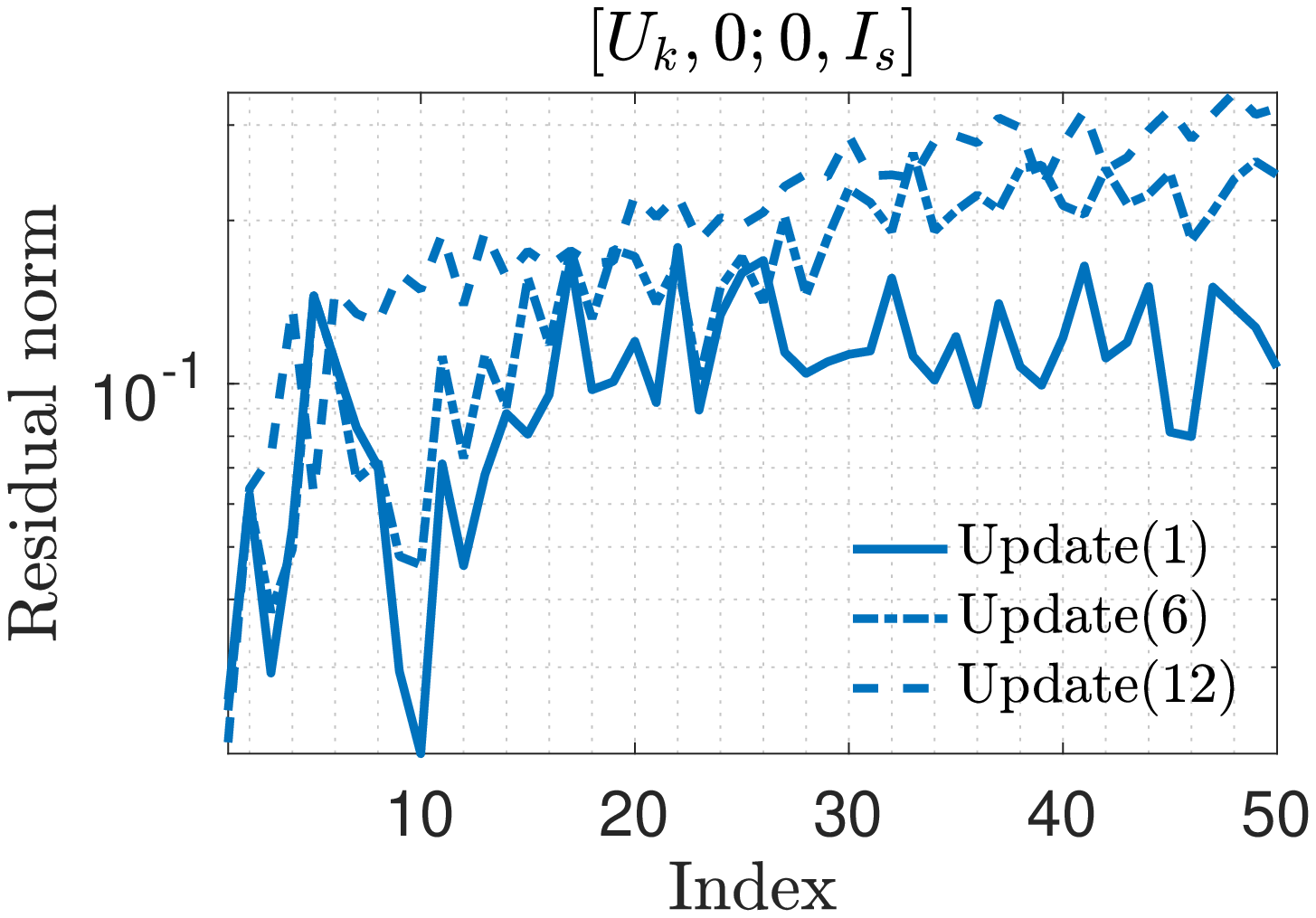}
     \includegraphics[width=0.24\linewidth]{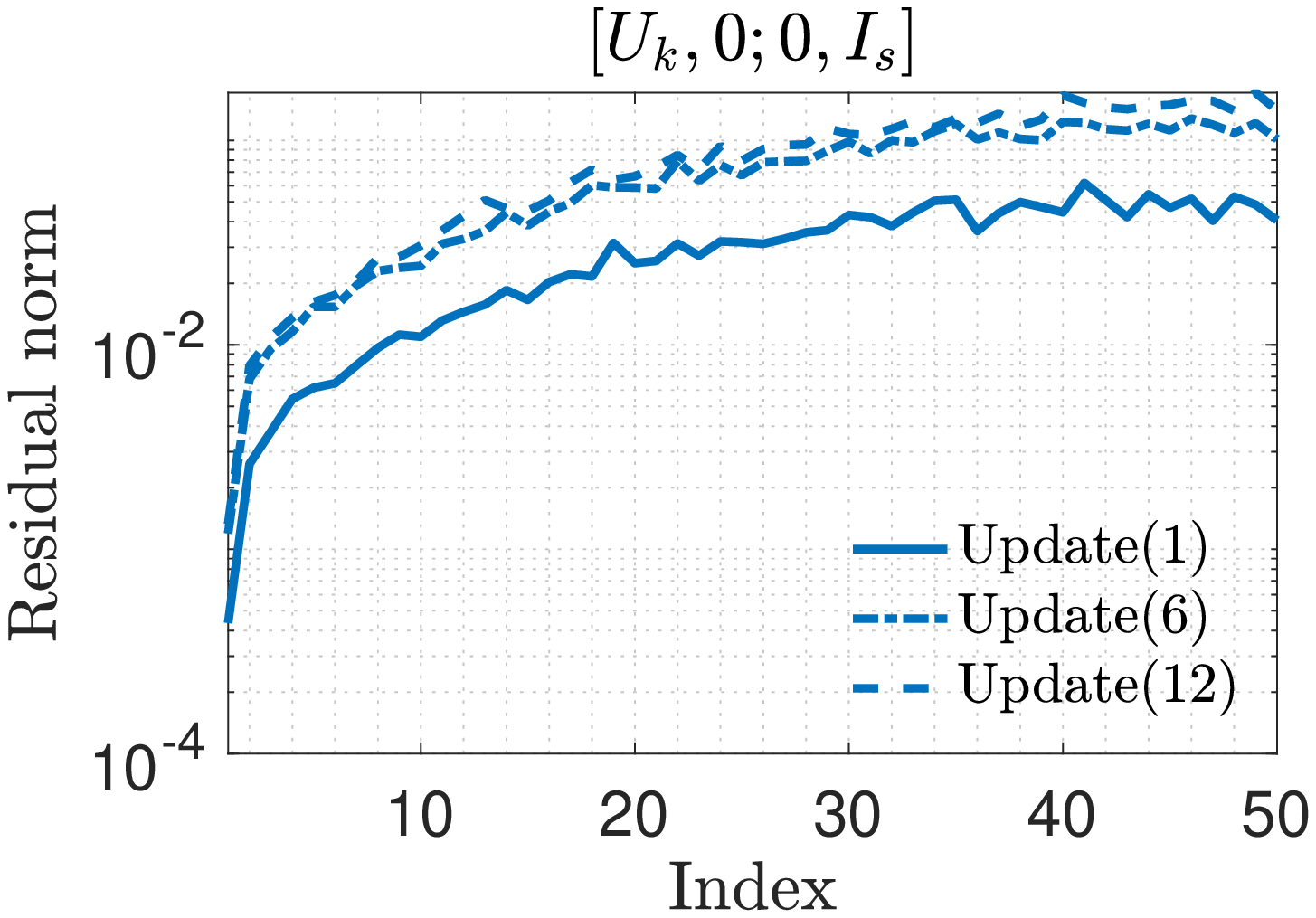}
     \includegraphics[width=0.24\linewidth]{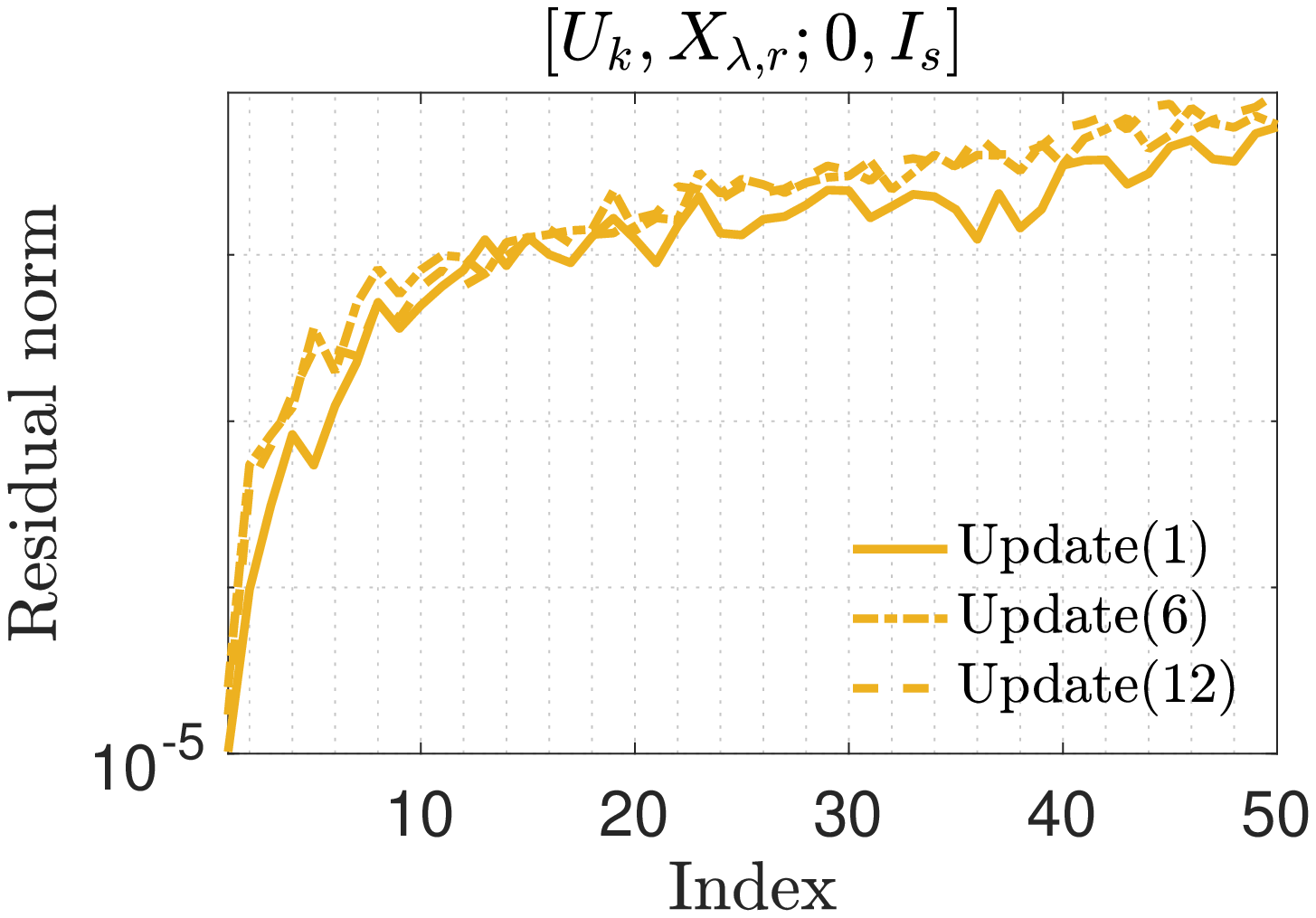}
     \includegraphics[width=0.24\linewidth]{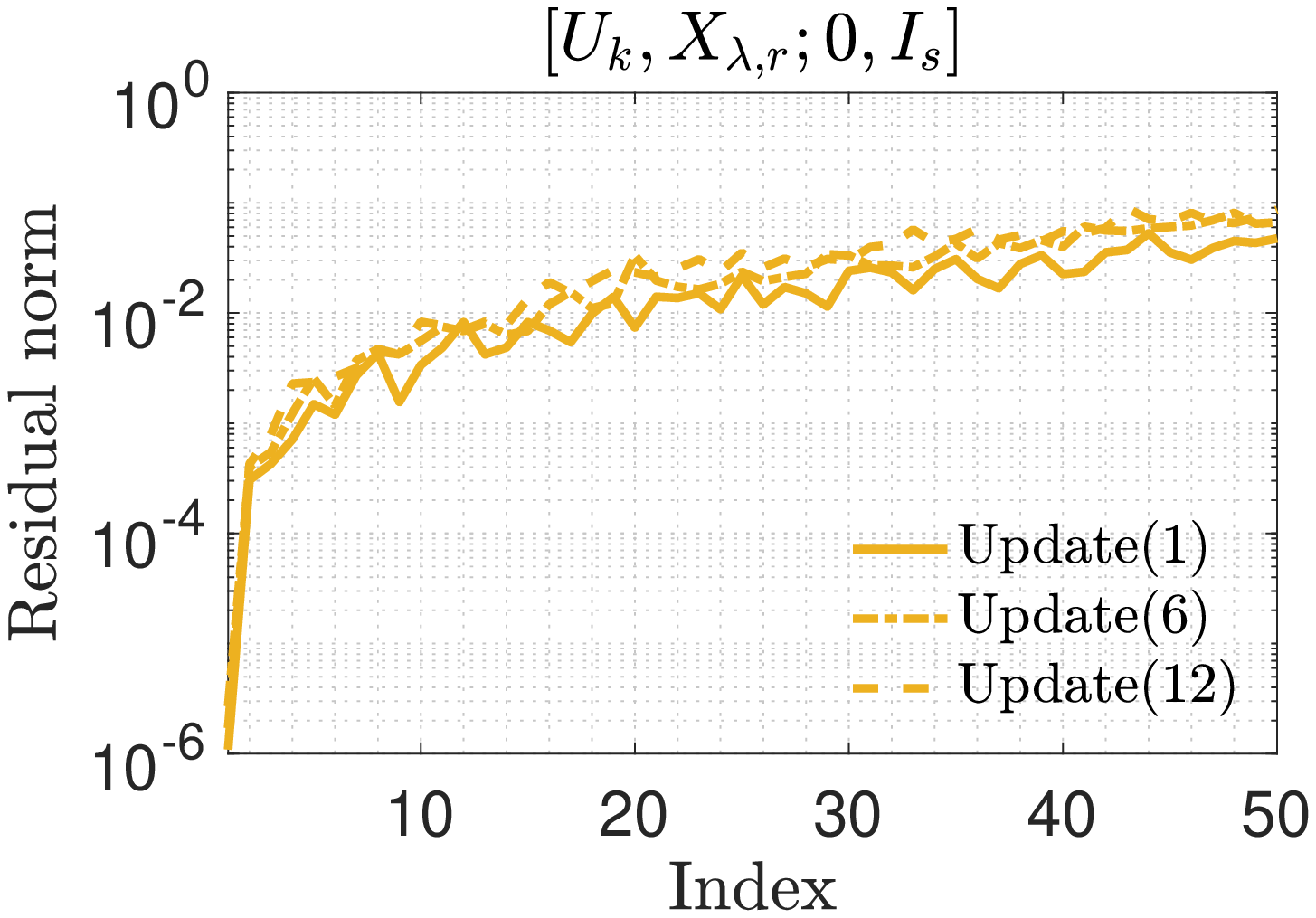}
     \includegraphics[width=0.24\linewidth]{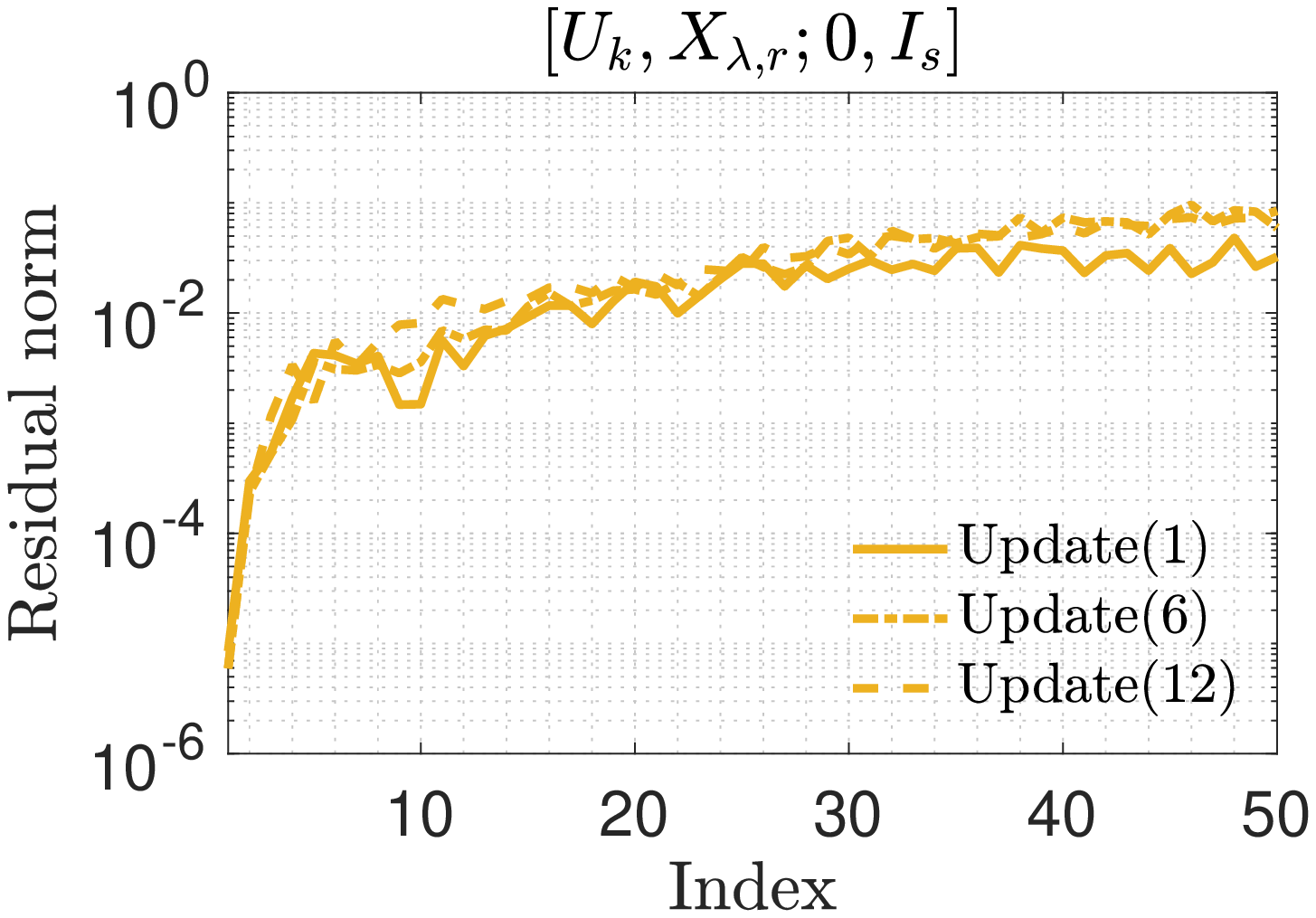}
     \includegraphics[width=0.24\linewidth]{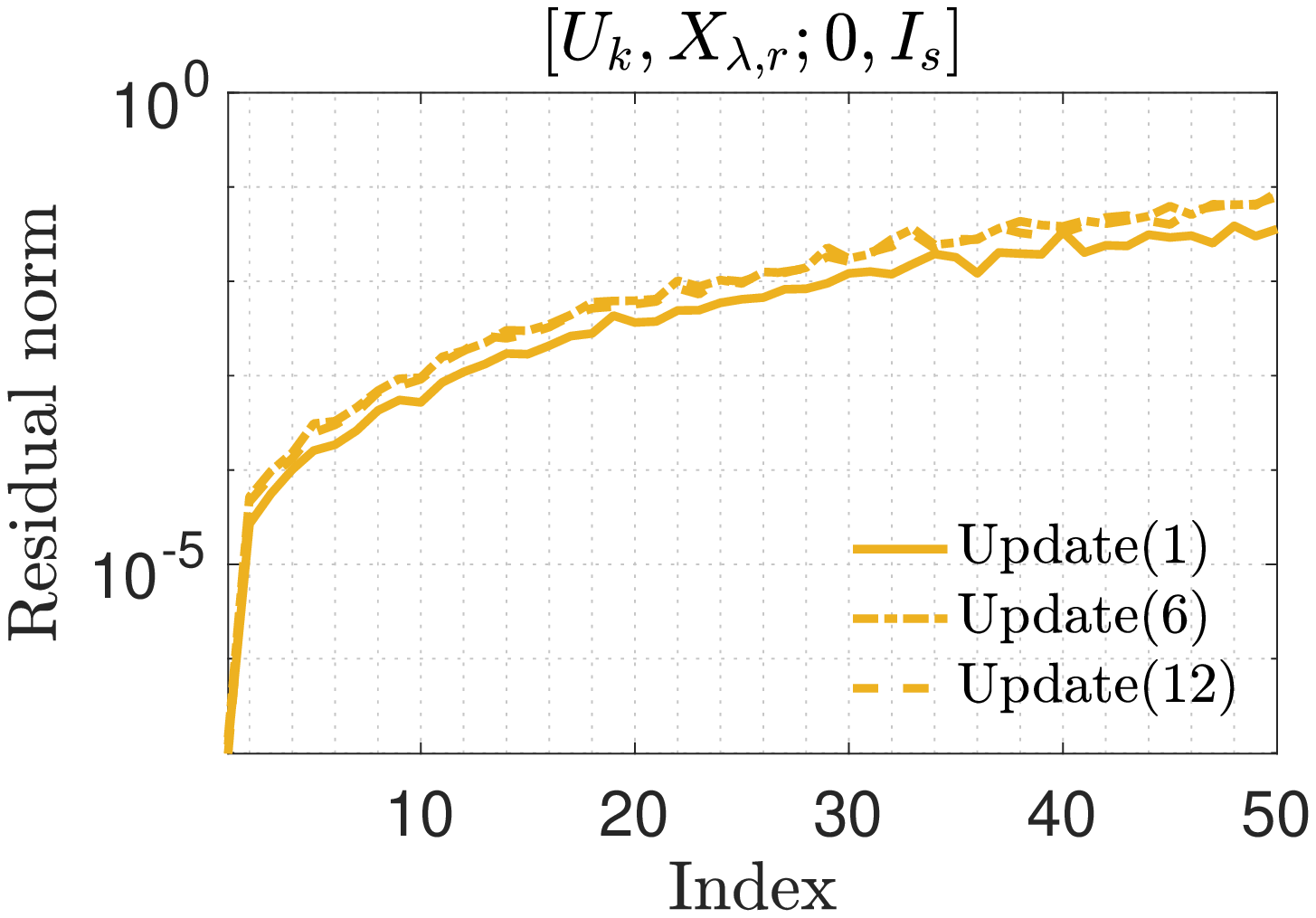}
     \caption{{\it Residual norm of the approximation of the $k=50$ leading 
     singular triplets of $A$ for the multiple updates case. From left to 
     right: MED, CRAN, CISI, and 1M.}}\label{fig:18}
\end{figure*}

\subsection{Sequence of updates.}

In this experiment the rows of matrix $E$ are now added in batches, i.e., 
we first approximate the $k$ leading singular triplets of matrix 
$A^{(0)} = \begin{pmatrix}
   B                                \\
   A(\ceil{m/2}+1:\ceil{m/2}+t,:)   \\
\end{pmatrix}$, then of matrix 
$A^{(1)} = \begin{pmatrix}
   A^{(0)}                              \\
   A(\ceil{m/2}+t+1:\ceil{m/2}+2t,:)    \\
\end{pmatrix}$, etc. Here, $t=\ceil{m/2}/\phi$ denotes the step-size and 
$\phi \in \mathbb{Z}^*$ denotes the total number of updates. Note that 
after the first update, the matrices $U_k$ and $V_k$ no longer denote 
the exact $k$ leading left and right singular vectors of the 
$B\equiv A^{(j-1)}$ submatrix of matrix $A^{(j)}$. We set $\phi=12$ and 
plot the accuracy achieved after one, six, and twelve updates, in Figures 
\ref{fig:17} and \ref{fig:18}. Notice that enhancing $Z$ by $X_{\lambda,r}$ 
leads to similar accuracy for all updates, while in the opposite case accuracy 
deteriorates as the updates accumulate. On a separate note, the accuracy 
of the $k$ leading singular triplets of $A$ is higher when matrix $E$ 
is added to $B$ in batches rather than in a single update as in the previous 
section.


\begin{table*}[!t]
\centering
\caption{\it Maximum relative error and residual norm of the approximation 
of the $k$ leading singular triplets of $A$ for the multiple updates case 
as $k$ varies. \label{table3}}
\vspace{0.05in}
\begin{tabular}{l @{\hskip 0.3in} r  c c  c  c c c c c c c c}
\toprule
\toprule
\multirow{2}{*}{} & 
& \multicolumn{2}{c}{\textbf{MED}} & &
\multicolumn{2}{c}{\textbf{CRAN}} & & \multicolumn{2}{c}{\textbf{CISI}} & & \multicolumn{2}{c}{\textbf{ML1M}}\\
\cmidrule[0.4pt](lr{0.125em}){3-4}%
\cmidrule[0.4pt](lr{0.125em}){6-7}%
\cmidrule[0.4pt](lr{0.125em}){9-10}%
\cmidrule[0.4pt](lr{0.125em}){12-13}%
&  Method & err. & res. & & err. & res. & & err. & res. & & err. & res. \\
\midrule
\multirow{2}{*}{$k=10$} & \cite{zha1999updating} & 0.046 & 0.172  & & 0.043 & 0.192  & & 0.054 & 0.274 &  & 0.002 &  0.058\\
 & \tikzmarkin[hor=style mygrey]{el10} Alg. \ref{alg1}                               & 0.001 & 0.045 &  & 0.008 & 0.090 & & 0.002 & 0.054 & &
 $3.0\mathtt{e}$-5 & 0.007 \tikzmarkend{el10} \\ \midrule
 \multirow{2}{*}{$k=20$} & \cite{zha1999updating} & 0.067 & 0.212 &  & 0.064 & 0.255 &  & 0.075 & 0.224 &   & 0.022 & 0.131 \\
 & \tikzmarkin[hor=style mygrey]{el20} Alg. \ref{alg1}                                & 0.004 & 0.073 &  & 0.005 & 0.076 &  & 0.003 & 0.053 &  & 0.002 & 0.040 \tikzmarkend{el20}\\ \midrule
 \multirow{2}{*}{$k=30$} & \cite{zha1999updating} & 0.076 & 0.384  & & 0.060 & 0.290  & & 0.084 & 0.330 &  & 0.023 & 0.123 \\
 & \tikzmarkin[hor=style mygrey]{el30} Alg. \ref{alg1}                                  & 0.006 & 0.067  & & 0.008 & 0.088  & & 0.004 & 0.070 &  & 0.001 & 0.041 \tikzmarkend{el30}\\ \bottomrule \bottomrule
\end{tabular}
\end{table*}

Table \ref{table3} lists relative error and residual norm associated 
with the approximation of the singular triplet $(\widehat{\sigma}_{50},\widehat{u}^{(50)},\widehat{v}^{(50)})$ 
by Algorithm \ref{alg1} and the method in \cite{zha1999updating}. 
The number of sought singular triplets $k$ was varied from ten to 
thirty. Comparisons against the method in \cite{vecharynski2014fast} 
were also performed but not reported since the latter was always less 
accurate than \cite{zha1999updating}. Overall, Algorithm \ref{alg1} 
provided higher accuracy, especially for those singular triplets whose 
corresponding singular value was closer to $\lambda$.

\section{Conclusion}

This paper presented an algorithm to update the rank-$k$ truncated SVD
of evolving matrices. The proposed algorithm undertakes a projection 
viewpoint and aims on building a pair of subspaces which approximate 
the linear span of the $k$ leading singular vectors of the updated 
matrix. Two different options to set these subspaces were considered. 
Experiments performed on matrices stemming from applications in LSI and 
recommender systems verified the effectiveness of the proposed scheme 
in terms of accuracy. 

\newpage 

\bibliographystyle{siam} 

\newpage
\onecolumn
\appendix

\begin{center}
\LARGE 
    Supplementary Material
\end{center}
\bigskip

\section*{\bf \Large Asymptotic complexity}

The asymptotic complexity analysis of the method in \cite{zha1999updating} is as follows. 
We need $O\left(ns^2 + nsk\right)$ FLOPs to form $(I_s-V_kV_k^H)E^H$ and compute its QR decomposition. The SVD of the matrix $Z^HAW$ requires $O\left((k+s)^3\right)$ FLOPs. 
Finally, the cost to form the approximation of matrices $\widehat{U}_k$ and $\widehat{V}_k$ 
is equal to $O\left(k^2(m+n) + nsk\right)$ FLOPs. 

The asymptotic complexity analysis for the “SV" variant of the method in \cite{vecharynski2014fast} 
is as follows. We need $O\left((\mathtt{nnz}(E)+nk)\delta_1+(n+s)\delta_1^2\right)$ FLOPs to approximate the $r$ leading singular triplets of $(I_s-V_kV_k^H)E^H$, where $\delta_1 \in \mathbb{Z}^*$ is greater than or equal to $r$ (i.e., $\delta_1$ is the number of Lanczos bidiagonalization steps). The cost to form and compute the SVD of the matrix $Z^HAW$ is equal to $(k+s)(k+r)^2 + \mathtt{nnz}(E)k+rs$ where the first term stands for the actual SVD and the rest of the terms stand for the formation of the matrix $Z^HAW$. Finally, the cost to form the approximation of matrices $\widehat{U}_k$ and $\widehat{V}_k$ is equal to $O\left(k^2(m+n) + nrk\right)$ FLOPs.

The asymptotic complexity analysis of Algorithm \ref{alg1} is as follows. 
First, notice that Algorithm \ref{alg1} requires no effort to build $W$.
For the case where $Z$ is set as in Proposition \ref{pro2}, termed as 
“Alg. \ref{alg1} (a)", we also need no FLOPs to build $Z$. The cost to 
solve the projected problem by unrestarted Lanczos is then equal to $O\left((\mathtt{nnz}(E)+nk)\delta_2+(k+s)\delta_2^2\right)$ FLOPs, 
where $\delta_2 \in \mathbb{Z}^*$ is greater than or equal to $k$ (i.e., 
$\delta_2$ is the number of unrestarted Lanczos steps). Finally, the cost 
to form the approximation of matrices $\widehat{U}_k$ and $\widehat{V}_k$ 
is equal to $O\left(k^2m+(\mathtt{nnz}(A)+n)k\right)$ FLOPs. For the case 
where $Z$ is set as in Proposition \ref{pro35}, termed as “Alg. \ref{alg1} 
(b)", we need 
\begin{equation*}
    \chi = O\left(\mathtt{nnz}(A)\delta_3 +m\delta_3^2\right)
\end{equation*}
FLOPs to build $X_{\lambda,r}$, where $\delta_3 \in \mathbb{Z}^*$ is greater 
than or equal to $k$ (i.e., $\delta_3$ is either the number of Lanczos 
bidiagonalization steps or the number of columns of matrix $R$ in randomized SVD). 

\begin{table}[ht]
\small 
\centering
\caption{\it Detailed asymptotic complexity of Algorithm \ref{alg1} and the schemes
in \cite{zha1999updating} and \cite{vecharynski2014fast}. All $\delta$ variables are replaced by $k$. \label{table5}}\vspace{0.01in}
\begin{tabular}{ l c c c c}
  \toprule
  \toprule
  Scheme & Building $Z$ & Building $W$ & Solving the projected problem & Other \\
  \midrule
  \myrowcolour
  \cite{zha1999updating}     & - & $ns^2 + nsk$ & $(k+s)^3$ & $k^2(m+n) + nsk$\\
  \cite{vecharynski2014fast} & - & $(\mathtt{nnz}(E)+nk)k+(n+s)k^2$ & $(k+s)(k+r)^2 + \mathtt{nnz}(E)k+rs$ & $k^2(m+n) + nrk$\\
  \myrowcolour
  Alg. \ref{alg1} (a)   & - & - & $(\mathtt{nnz}(E)+nk)k+(k+s)k^2$ & $k^2m+(\mathtt{nnz}(A)+n)k$ \\
  Alg. \ref{alg1} (b)   & $\chi$ & - & $(\mathtt{nnz}(E)+(n+r)k)k+(k+r+s)k^2$ & $k^2m+(\mathtt{nnz}(A)+n)k$ \\
  \bottomrule
  \bottomrule
\end{tabular}
\end{table}

The above discussion is summarized in Table \ref{table5} where we list the 
asymptotic complexity of Algorithm \ref{alg1} and the schemes
in \cite{zha1999updating} and \cite{vecharynski2014fast}. The complexities 
of the latter two schemes were also verified by adjusting the complexity 
analysis from \cite{vecharynski2014fast}. To allow for a practical comparison, 
we replaced all $\delta$ variables with $k$ since in practice these variables 
are equal to at most a small integer multiple of $k$.

Consider now a comparison between Algorithm \ref{alg1} (a) and the method in  
\cite{zha1999updating}. For all practical purposes, these two schemes return 
identical approximations to $A_k$. Nonetheless, Algorithm \ref{alg1} (a) 
requires no effort to build $W$. Moreover, the cost to solve the projected 
problem is linear with respect to $s$ and cubic with respect to $k$, instead 
of cubic with respect to the sum $s+k$ in \cite{zha1999updating}. The only 
scenario where Algorithm \ref{alg1} can be potentially more expensive than 
\cite{zha1999updating} is when matrix $A$ is exceptionally dense, and both 
$k$ and $s$ are very small. Similar observations can be made for the relation 
between Algorithm \ref{alg1} (b) and the methods in \cite{vecharynski2014fast}, 
although the comparison is more involved.

\section*{\bf \Large Proofs}

\subsection*{Proof of Proposition \ref{pro0}}

The scalar-vector pair $(\widehat{\sigma}_{i}^2,\widehat{u}^{(i)})$ satisfies the 
equation $(AA^H-\widehat{\sigma}_{i}^2 I_{m+s})\widehat{u}^{(i)}=0$. If we partition 
the $i$'th left singular vector as 
\[
\widehat{u}^{(i)} =
 \begin{pmatrix}
    \widehat{f}^{(i)}          \\[0.3em]
    \widehat{y}^{(i)}          \\[0.3em]
  \end{pmatrix},
\]
we can write
\begin{equation*} 
  \begin{pmatrix}
   BB^H-\widehat{\sigma}_{i}^2 I_{m} & BE^H            \\[0.3em]
   EB^H & EE^H-\widehat{\sigma}_{i}^2 I_{s}            \\[0.3em]
  \end{pmatrix}
  \begin{pmatrix}
   \widehat{f}^{(i)}        \\[0.3em]
   \widehat{y}^{(i)}        \\[0.3em]
  \end{pmatrix}=0.
\end{equation*}
The leading $m$ rows satisfy $(BB^H-\widehat{\sigma}_{i}^2I_m)\widehat{f}^{(i)}=-BE^H\widehat{y}^{(i)}$. 
Plugging the expression of $\widehat{f}^{(i)}$ in the second block of rows and considering the 
full SVD $B=U \Sigma V^H$ leads to 
\begin{align*}
0 & = \left[EE^H -EB^H(BB^H-\widehat{\sigma}_{i}^2I_m)^{-1}BE^H-\widehat{\sigma}_{i}^2I_s\right]
\widehat{y}^{(i)} \\
  & =  \left[E(I_s-B^H(BB^H-\widehat{\sigma}_{i}^2I_m)^{-1}B)E^H-\widehat{\sigma}_{i}^2I_s\right]
  \widehat{y}^{(i)}\\
  & =  \left[E(VV^H+V\Sigma^T(\widehat{\sigma}_i^2 I_m - \Sigma \Sigma^T)^{-1} \Sigma V^H)E^H-\widehat{\sigma}_{i}^2I_s\right]\widehat{y}^{(i)}\\
  & =  \left[EV(I_n+\Sigma^T\left(\widehat{\sigma}_i^2 I_m - \Sigma \Sigma^T\right)^{-1} \Sigma)V^HE^H-\widehat{\sigma}_{i}^2I_s\right]\widehat{y}^{(i)}.
\end{align*}
The proof concludes by noticing that 
\begin{equation*}
    I_n+\Sigma^T\left(\widehat{\sigma}_i^2 I_m - \Sigma \Sigma^T\right)^{-1} \Sigma
    =
    \begin{pmatrix}
    1+\dfrac{\sigma_1^2}{\widehat{\sigma}_i^2-\sigma_1^2} & & \\[0.3em]
    & \ddots & \\[0.3em]
    & & 1+\dfrac{\sigma_{n}^2}{\widehat{\sigma}_i^2-\sigma_{n}^2} \\[0.3em]
    \end{pmatrix}
    =
    \begin{pmatrix}
    \dfrac{\widehat{\sigma}_i^2}{\widehat{\sigma}_i^2-\sigma_1^2} & & \\[0.3em]
    & \ddots & \\[0.3em]
    & & \dfrac{\widehat{\sigma}_{i}^2}{\widehat{\sigma}_i^2
    -\sigma_{n}^2} \\[0.3em]
    \end{pmatrix},
\end{equation*}
where for the case $m<n$, we have $\sigma_j=0$ for any $j=m+1,\ldots,n$. 
In case $\widehat{\sigma}_{i}=\sigma_{j}$, the Moore-Penrose 
pseudoinverse $(BB^H-\widehat{\sigma}_{i}^2I_m)^\dagger$ is 
considered instead.

\subsection*{Proof of Proposition \ref{pro1}}
Since the left singular vectors of $B$ span $\mathbb{R}^m$, we can write 
\begin{equation*}
    BE^H\widehat{y}^{(i)} = \sum\limits_{j=1}^m \sigma_{j} 
    u^{(j)} \left(Ev^{(j)}\right)^H\widehat{y}^{(i)}. 
\end{equation*}
The proof concludes by noticing that the top $m\times 1$ part of 
$\widehat{u}^{(i)}$ can be written as 
    \begin{align*}
        \widehat{f}^{(i)} & = -(BB^H-\widehat{\sigma}_{i}^2I_m)^{-1}BE^H\widehat{y}^{(i)} \\
        &= -U(\Sigma \Sigma^T - \widehat{\sigma}_i^2 I_m)^{-1}\Sigma \left(EV\right)^H\widehat{y}^{(i)}\\
        &=-\sum\limits_{j=1}^{\mathtt{min}(m,n)} u^{(j)}\dfrac{\sigma_{j}}{\sigma_{j}^2-\widehat{\sigma}_{i}^2}  
        \left(Ev^{(j)}\right)^H\widehat{y}^{(i)}\\
        &=-\sum\limits_{j=1}^{\mathtt{min}(m,n)}u^{(j)}
        \dfrac{\sigma_{j}}{\sigma_{j}^2-\widehat{\sigma}_{i}^2} \left(Ev^{(j)}\right)^H\widehat{y}^{(i)}\\
        &=\sum\limits_{j=1}^{\mathtt{min}(m,n)}u^{(j)} \chi_{j,i}.
    \end{align*}

\subsection*{Proof of Proposition \ref{pro2}}
We have
        \begin{align*}
            \mathtt{min}_{z\in \mathtt{range}(Z)} \|\widehat{u}^{(i)}-z\| 
            & \leq 
            \norm{   \begin{pmatrix}
             u^{(k+1)},\ldots,u^{(\mathtt{min}(m,n))} \\[0.3em]
                            \\[0.3em]
  \end{pmatrix}
  \begin{pmatrix}
   \chi_{k+1,i}   \\[0.3em]
   \vdots         \\[0.3em]
   \chi_{\mathtt{min}(m,n),i}     \\[0.3em]
  \end{pmatrix}} \\
            & = 
  \norm{\left( 
  \begin{smallmatrix}
  \scalebox{2}{$0$}_{k,k}  & & & & \\ 
  &   & \dfrac{\sigma_{k+1}}{\sigma_{k+1}^2-\widehat{\sigma}_{i}^2} & & \\ 
  &   & & \ddots & \\ 
  &  & & & \dfrac{\sigma_{\mathtt{min}(m,n)}}{\sigma_{\mathtt{min}(m,n)}^2-\widehat{\sigma}_{i}^2}\\
 \end{smallmatrix}
 \right)V^HE^H\hat{y}^{(i)}}
  \\
            & \leq  
\mathtt{max}\left\lbrace\left|\dfrac{\sigma_{j}}{\sigma_{j}^2-\widehat{\sigma}_{i}^2}\right|\right\rbrace_{j=k+1,\ldots,\mathtt{min}(m,n)}
            \norm{E^H\widehat{y}^{(i)}}.
        \end{align*}
    The proof follows by noticing that due to Cauchy's interlacing theorem 
    we have $\sigma_{k+1}^2\leq \widehat{\sigma}_{i}^2,\ i=1,\ldots,k$, and 
    thus $\left|\dfrac{\sigma_{k+1}}{\sigma_{k+1}^2-\widehat{\sigma}_{i}^2}\right|
    \geq \cdots \geq \left|\dfrac{\sigma_{\mathtt{min}(m,n)}}
    {\sigma_{\mathtt{min}(m,n)}^2-\widehat{\sigma}_{i}^2}\right|.$

\subsection*{Proof of Lemma \ref{lem1}}

We can write 
\begin{equation*}
\begin{aligned}
B(\lambda) &=\left(I-U_k U_k^H\right) U
\left(\begin{smallmatrix} \sigma_1^2-\lambda & &\\ & \ddots &\\ & & \sigma_m^2-\lambda \end{smallmatrix}\right)^{-1}U^H\\
 &= 
  U
  \left( 
  \begin{smallmatrix}
  \scalebox{2}{$0$}_{k,k}  & & & & \\ 
  &   & \dfrac{1}{\sigma_{k+1}^2-\lambda} & & \\ 
  &   & & \ddots & \\ 
  &  & & & \dfrac{1}{\sigma_m^2-\lambda}\\
 \end{smallmatrix}
 \right)U^{H},
 \end{aligned}
 \end{equation*} 
 where $\sigma_j=0$ for any $j>\mathtt{min}(m,n)$. Let us now define the scalar 
$\gamma_{j,i} = \dfrac{\widehat{\sigma}_i^2-\lambda}{\sigma_j^2-\lambda}$. Then,

\begin{equation*}
    B(\lambda)\left[(\widehat{\sigma}_i^2-\lambda)B(\lambda)\right]^\rho = U
    \begin{pmatrix}
    \scalebox{2}{$0$}_{k,k} & & & \\[0.3em]
    & \dfrac{\gamma_{k+1,i}^\rho}{\sigma_{k+1}^2-\lambda} & & \\[0.3em]
    & & \ddots & \\[0.3em]
    & & & \dfrac{\gamma_{m,i}^\rho}{\sigma_{m}^2-\lambda} \\[0.3em]
    \end{pmatrix}
    U^{H}.
\end{equation*}
Accounting for all powers $p=0,1,2,\ldots$, gives
{\small \begin{equation*}
    B(\lambda)\sum_{\rho=0}^{\infty} 
    \left[(\widehat{\sigma}_i^2-\lambda)B(\lambda)\right]^\rho = U
    \begin{pmatrix}
    \scalebox{2}{$0$}_{k,k} & & & \\[0.3em]
    & \dfrac{\sum_{\rho=0}^{\infty} \gamma_{k+1,i}^\rho}{\sigma_{k+1}^2-\lambda} & & \\[0.3em]
    & & \ddots & \\[0.3em]
    & & & \dfrac{\sum_{\rho=0}^{\infty} \gamma_{m,i}^\rho}{\sigma_{m}^2-\lambda} \\[0.3em]
    \end{pmatrix}
    U^{H}.
\end{equation*}}
Since $\lambda > \widehat{\sigma}_k^2 \geq \sigma_k^2$, it follows 
that for any $j>k$ we have $|\gamma_{j,i}| < 1$. 
Therefore, the geometric series converges and $\sum_{\rho=0}^{\infty} \gamma_{j,i}^\rho=\dfrac{1}{1-\gamma_{j,i}}=
\dfrac{\sigma_j^2-\lambda}{\sigma_j^2-\widehat{\sigma}_i^2}$. 
It follows that $\dfrac{1}{\sigma_j^2-\lambda}\sum_{\rho=0}^{\infty} \gamma_j^\rho
=\dfrac{1}{\sigma_j^2-\widehat{\sigma}_i^2}$.

We finally have
\begin{equation*}
\begin{aligned}
    B(\lambda)\sum_{\rho=0}^{\infty} 
    \left[(\widehat{\sigma}_i^2-\lambda)B(\lambda)\right]^\rho & = U
    \begin{pmatrix}
    \scalebox{2}{$0$}_{k,k} & & & & \\[0.3em]
    & \dfrac{1}{\sigma_{k+1}^2-\widehat{\sigma}_i^2} & & \\[0.3em]
    & & \ddots & \\[0.3em]
    & & & \dfrac{1}{\sigma_{m}^2-\widehat{\sigma}_i^2} \\[0.3em]
    \end{pmatrix}
    U^{H} \\
    & = \left(I-U_k U_k^H\right)B(\widehat{\sigma}_i^2).
    \end{aligned}
\end{equation*}
This concludes the proof.

\subsection*{Proof of Proposition \ref{pro34}} \label{proof1}

    First, notice that 
    \begin{equation*}
    (BB^H-\widehat{\sigma}_i^2 I_m)^{-1} = 
    U_k U_k^H(BB^H-\widehat{\sigma}_i^2 I_m)^{-1} + (I_m -U_k U_k^H) (BB^H-\widehat{\sigma}_i^2 I_m)^{-1}.
    \end{equation*}
    Therefore, we can write 
    \[(BB^H-\widehat{\sigma}_i^2 I_m)^{-1}BE^H\widehat{y}^{(i)} 
    = U_k(\Sigma_k^2-\widehat{\sigma}_i^2 I_k)^{-1} \Sigma_k (EV_k)^H\widehat{y}^{(i)}
    +(I_m -U_k U_k^H) (BB^H-\widehat{\sigma}_i^2 I_m)^{-1}BE^H\widehat{y}^{(i)}.\]
    The left singular vector $\widehat{u}^{(i)}$ can be then expressed as
   \begin{eqnarray*}
   \widehat{u}^{(i)} & = &
   \begin{pmatrix}
   -(BB^H-\widehat{\sigma}_i^2 I_m)^{-1}BE^H                         \\[0.3em]
   I_s                      \\[0.3em]
   \end{pmatrix}\widehat{y}^{(i)} \\
   & = &
   \begin{pmatrix}
   u^{(1)},\ldots,u^{(k)} & \\[0.3em]
   & I_s                    \\[0.3em]
   \end{pmatrix}
   \begin{pmatrix}
   \chi_{1,i}   \\[0.3em]
   \vdots       \\[0.3em]
   \chi_{k,i}   \\[0.3em]
   \widehat{y}^{(i)}      \\[0.3em]
   \end{pmatrix} - 
   \begin{pmatrix}
    B(\widehat{\sigma}_i^2)BE^H\widehat{y}^{(i)}  \\[0.3em]
                                            \\[0.3em]
  \end{pmatrix}.
  \end{eqnarray*}
  The proof concludes by noticing that by Lemma \ref{lem1} we have $B(\widehat{\sigma}_i^2)
  =B(\lambda) \sum\limits_{\rho=0}^\infty 
   \left[(\widehat{\sigma}_i^2-\lambda)B(\lambda)\right]^\rho$. 

\subsection*{Proof of Proposition \ref{pro35}}

The proof exploits the formula  
\begin{equation*}
    (B(\widehat{\sigma}_i^2)-B(\lambda))BE^H = 
    (I-U_kU_k^H)U\left[(\Sigma \Sigma^T-\widehat{\sigma}_i^2 I_m)^{-1}-(\Sigma \Sigma^T-\lambda I_m)^{-1}\right]U^HU\Sigma V^HE^H.
\end{equation*}
It follows
    \begin{eqnarray*}
            \mathtt{min}_{z\in \mathtt{range}(Z)} \|\widehat{u}^{(i)}-z\| 
            & \leq  & 
            \norm{\begin{pmatrix}
            \left[B(\widehat{\sigma}_i^2)-B(\lambda)\right]BE^H\widehat{y}^{(i)}\\[0.3em]
                                                                         \\[0.3em]
            \end{pmatrix}}\\
            & \leq &
            \norm{   \begin{pmatrix}
    \scalebox{2}{$0$}_{k,k} & & & & \\[0.3em]
    & \dfrac{\sigma_{k+1}(\widehat{\sigma}_i^2-\lambda)}
            {(\sigma_{k+1}^2-\widehat{\sigma}_{i}^2)\left(\sigma_{k+1}^2-\lambda\right)} & & \\[0.3em]
    & & \ddots & \\[0.3em]
    & & & \dfrac{\sigma_{\mathtt{min}(m,n)}(\widehat{\sigma}_i^2-\lambda)}
            {(\sigma_{\mathtt{min}(m,n)}^2-\widehat{\sigma}_{i}^2)\left(\sigma_{\mathtt{min}(m,n)}^2-\lambda\right)} \\[0.3em]
    \end{pmatrix}}\norm{E^H\widehat{y}^{(i)}}\\
            & \leq &
          \mathtt{max}\left\lbrace\left|\dfrac{\sigma_{j}(\widehat{\sigma}_i^2-\lambda)}
            {(\sigma_{j}^2-\widehat{\sigma}_{i}^2)\left(\sigma_{j}^2-\lambda\right)}
            \right|\right\rbrace_{j=k+1,\ldots,\mathtt{min}(m,n)}
            \norm{E^H\widehat{y}^{(i)}}. 
        \end{eqnarray*}

\end{document}